# MAXIMUM EMPIRICAL LIKELIHOOD ESTIMATION OF THE SPECTRAL MEASURE OF AN EXTREME-VALUE DISTRIBUTION

By John H. J. Einmahl and Johan Segers[1]

*Tilburg University and Université catholique de Louvain*

Consider a random sample from a bivariate distribution function $F$ in the max-domain of attraction of an extreme-value distribution function $G$. This $G$ is characterized by two extreme-value indices and a spectral measure, the latter determining the tail dependence structure of $F$. A major issue in multivariate extreme-value theory is the estimation of the spectral measure $\Phi_p$ with respect to the $L_p$ norm. For every $p \in [1, \infty]$, a nonparametric maximum empirical likelihood estimator is proposed for $\Phi_p$. The main novelty is that these estimators are guaranteed to satisfy the moment constraints by which spectral measures are characterized. Asymptotic normality of the estimators is proved under conditions that allow for tail independence. Moreover, the conditions are easily verifiable as we demonstrate through a number of theoretical examples. A simulation study shows a substantially improved performance of the new estimators. Two case studies illustrate how to implement the methods in practice.

**1. Introduction.** Let $F$ be a continuous bivariate distribution function in the max-domain of attraction of an extreme-value distribution function $G$. Up to location and scale, the marginals of $G$ are determined by the extreme-value indices of the marginals of $F$. The dependence structure of $G$ can be described in various equivalent ways; in this paper we focus on the spectral measure $\Phi$ introduced in de Haan and Resnick (1977). The spectral

[1]Research Fellow of CentER, Tilburg University. Supported by the Netherlands Organization for Scientific Research (NWO) in the form of a VENI grant and by the IAP research network Grant P6/03 of the Belgian government (Belgian Science Policy).

or angular measure is a finite Borel measure on a compact interval, here taken to be $[0, \pi/2]$.

Given a random sample from $F$, statistical inference on the upper tail of $F$ falls apart into two pieces: estimation of the upper tails of its marginal distributions, which is well understood, and estimation of $\Phi$, which we will consider in this paper. The actual representation of the spectral measure depends on the norm used on $\mathbb{R}^2$; here we will consider the $L_p$ norm for every $p \in [1, \infty]$, with $\Phi_p$ denoting the corresponding spectral measure. The most common choices in the literature are $p = 1$, 2 and $\infty$.

It is the aim of this paper to derive a nonparametric estimator of the spectral measure, superior to its predecessors, and to establish its asymptotic normality. In Einmahl, de Haan and Piterbarg (2001), a nonparametric estimator $\hat{\Phi}_\infty$ was proposed for $\Phi_\infty$. This estimator, which we will refer to as the empirical spectral measure, was shown to be asymptotically normal under the assumption that $\Phi_\infty$ has a density, excluding thereby the case of tail independence. Moreover the empirical spectral measure is itself not a proper spectral measure because it violates the moment constraints characterizing the class of spectral measures. A related estimator in a more restrictive framework was proposed in Einmahl, de Haan and Sinha (1997).

The contributions of our paper are threefold: first, to propose a nonparametric estimator for the spectral measure which itself satisfies the moment constraints; second, to allow for arbitrary $L_p$ norms, $p \in [1, \infty]$; third, to prove asymptotic normality under flexible and easily verifiable conditions that allow for spectral measures with atoms at 0 or $\pi/2$, including thereby the case of tail independence. We do this in two steps: first we define for every $p \in [1, \infty]$ the empirical spectral measure $\hat{\Phi}_p$ and extend the results in Einmahl, de Haan and Piterbarg (2001) under the weaker conditions mentioned above; second, we use a nonparametric maximum empirical likelihood approach to enforce the moment constraints, thereby obtaining an estimator $\tilde{\Phi}_p$ that is itself a genuine spectral measure. A simulation study shows that the new estimator $\tilde{\Phi}_p$ is substantially more efficient than the empirical spectral measure $\hat{\Phi}_p$.

As the new estimator takes values in the class of spectral measures, it can be easily transformed into estimators for the aforementioned other objects that can be used to describe the dependence structure of $G$. This holds in particular for the Pickands (1981) dependence function and the stable tail dependence function [Huang (1992), Drees and Huang (1998) and Einmahl, de Haan and Li (2006)]. For a general background on spectral measures and these dependence functions as well as results for the corresponding estimators, see for instance the monographs Coles (2001), Beirlant et al. (2004) and de Haan and Ferreira (2006).

An alternative to the nonparametric approach in this paper is the parametric one [Coles and Tawn (1991) and Joe, Smith and Weissman (1992)].

Parametric models for the spectral measure are usually defined for $p = 1$ because this choice tends to lead to simpler formulae. Many parametric models, such as the asymmetric (negative) logistic and the asymmetric mixed models, allow the spectral measure to have atoms at 0 and $\pi/2$.

The paper is organized as follows. In Section 2 we review the general probabilistic theory for spectral measures. The asymptotic normality results for $\hat{\Phi}_p$ and $\tilde{\Phi}_p$ are presented in Sections 3 and 4, respectively. In Section 5 some theoretical examples are discussed and used in a simulation study; moreover, the methods are applied to a sample of insurance indemnity claims taken from Frees and Valdez (1998) and to body measurement data of the National Health and Nutrition Examination Survey 2005–2006 [National Center for Health Statistics (2007)]. Sections 6 and 7 contain the proofs of the results in Sections 3 and 4, respectively.

**2. Spectral measures.** Let $(X_1, X_2)$ be a bivariate random vector with continuous distribution function $F$ and marginal distribution functions $F_1$ and $F_2$. Put

$$Z_j = \frac{1}{1 - F_j(X_j)}, \qquad j = 1, 2. \tag{2.1}$$

Define $\mathbb{E} = [0, \infty]^2 \setminus \{(0, 0)\}$. Assume that

$$s \Pr[s^{-1}(Z_1, Z_2) \in \cdot] \xrightarrow{v} \mu(\cdot), \qquad s \to \infty, \tag{2.2}$$

where "$\xrightarrow{v}$" stands for *vague convergence* of measures (in $\mathbb{E}$): for every continuous $f : \mathbb{E} \to \mathbb{R}$ with compact support, $\lim_{s\to\infty} sE[f(s^{-1}(Z_1, Z_2))] = \int_{\mathbb{E}} f\, d\mu$.

The *exponent measure* $\mu$ enjoys two crucial properties: homogeneity,

$$\mu(c\cdot) = c^{-1}\mu(\cdot), \qquad 0 < c < \infty, \tag{2.3}$$

and standardized marginals,

$$\mu([z, \infty] \times [0, \infty]) = \mu([0, \infty] \times [z, \infty]) = 1/z, \qquad 0 < z \le \infty. \tag{2.4}$$

Note that $\mu$ is concentrated on $[0, \infty)^2 \setminus \{(0, 0)\}$, that is, $\mu([0, \infty]^2 \setminus [0, \infty)^2) = 0$.

Let $\|\cdot\|$ be an arbitrary norm on $\mathbb{R}^2$; for convenience, assume that $\|(1, 0)\| = 1 = \|(0, 1)\|$. Consider the following polar coordinates, $(r, \theta)$, of $(z_1, z_2) \in [0, \infty)^2 \setminus \{(0, 0)\}$:

$$\begin{aligned} r &= \|(z_1, z_2)\| \in (0, \infty), \\ \theta &= \arctan(z_1/z_2) \in [0, \pi/2]. \end{aligned} \tag{2.5}$$

As we will see later, the choice of radial coordinate $r$ through the norm has important implications; the choice of the angular coordinate $\theta$ is unimportant, that is, we could just as well have used $z_1/(z_1 + z_2) \in [0, 1]$ or $z_1/\|(z_1, z_2)\|$.

Given the exponent measure $\mu$ and using polar coordinates $(r,\theta)$ as in (2.5), define a Borel measure $\Phi$ on $[0,\pi/2]$ by

$$\Phi(\cdot) = \mu(\{(z_1,z_2) \in [0,\infty)^2 : r \geq 1, \theta \in \cdot\}). \tag{2.6}$$

The *spectral measure* $\Phi$ admits the following interpretation in terms of $(Z_1,Z_2)$ in (2.1):

$$s \Pr[\|(Z_1,Z_2)\| \geq s, \arctan(Z_1/Z_2) \in \cdot] \xrightarrow{v} \Phi(\cdot), \qquad s \to \infty. \tag{2.7}$$

By homogeneity of $\mu$, see (2.3), for every $\mu$-integrable $f : \mathbb{E} \to \mathbb{R}$,

$$\int_{\mathbb{E}} f \, d\mu = \int_{[0,\pi/2]} \int_0^\infty f(z_1(r,\theta), z_2(r,\theta)) r^{-2} \, dr \, \Phi(d\theta), \tag{2.8}$$

where $z_1(r,\theta) = r\sin\theta / \|(\sin\theta,\cos\theta)\|$ and $z_2(r,\theta) = r\cos\theta/\|(\sin\theta,\cos\theta)\|$ form the inverse of the polar transformation (2.5). By (2.8), in the polar coordinate system $(r,\theta)$, the exponent measure $\mu$ is a product measure $r^{-2}\,dr\,\Phi(d\theta)$. In particular, the exponent measure $\mu$ is completely determined by its spectral measure $\Phi$. The standardization constraints (2.4) on $\mu$ translate into moment constraints on $\Phi$:

$$\int_{[0,\pi/2]} \frac{\sin\theta}{\|(\sin\theta,\cos\theta)\|} \Phi(d\theta) = 1 = \int_{[0,\pi/2]} \frac{\cos\theta}{\|(\sin\theta,\cos\theta)\|} \Phi(d\theta). \tag{2.9}$$

Note that $X_1$ and $X_2$ are tail independent, that is, $s \Pr[Z_1 \geq s, Z_2 \geq s] \to 0$ as $s \to \infty$, if and only if $\mu$ is concentrated on the coordinate axes, or, equivalently, $\Phi$ is concentrated on $\{0, \pi/2\}$; in that case, $\Phi(\{0\}) = 1 = \Phi(\{\pi/2\})$. On the other hand, the variables $X_1$ and $X_2$ are completely tail dependent, that is, $s\Pr[Z_1 \geq s,\ Z_2 \geq s] \to 1$ as $s \to \infty$, if and only if $\mu$ is concentrated on the main diagonal, or, equivalently, $\Phi$ is concentrated on $\{\pi/4\}$; in that case, $\Phi(\{\pi/4\}) = \|(1,1)\|$. In general, the total mass $\Phi([0,\pi/2])$ of a spectral measure is finite but even for a fixed norm it can vary for different exponent measures $\mu$, with one exception: in case of the $L_1$ norm, by addition of the two constraints in (2.9), $\Phi([0,\pi/2]) = 2$ for every exponent measure $\mu$. The spectral measure was introduced in de Haan and Resnick (1977); for more details on the results in this section see Beirlant et al. (2004) and de Haan and Ferreira (2006).

Dividing the spectral measure $\Phi$ by its total mass yields a probability measure $Q$ on $[0,\pi/2]$

$$Q(\cdot) = \Phi(\cdot)/\Phi([0,\pi/2]), \tag{2.10}$$

which we coin the *spectral probability measure*. By (2.7)

$$\Pr[\arctan(Z_1/Z_2) \in \cdot \mid \|(Z_1,Z_2)\| \geq s] \xrightarrow{d} Q(\cdot), \qquad s \to \infty. \tag{2.11}$$

In words, $Q$ is the limit distribution of the angle $\arctan(Z_1/Z_2)$ when the radius $\|(Z_1, Z_2)\|$ is large. The moment constraints (2.9) on $\Phi$ are equivalent to the following moment constraint on $Q$:

$$\int_{[0,\pi/2]} \frac{\sin\theta}{\|(\sin\theta, \cos\theta)\|} Q(d\theta) = \int_{[0,\pi/2]} \frac{\cos\theta}{\|(\sin\theta, \cos\theta)\|} Q(d\theta) =: m(Q). \tag{2.12}$$

Conversely, we can reconstruct $\Phi$ from $Q$ by

$$\Phi(\cdot) = Q(\cdot)/m(Q). \tag{2.13}$$

The spectral *probability* measure $Q$ allows nonparametric maximum likelihood estimation; see Section 4. The estimator of $\Phi$ then follows through (2.13).

In Einmahl, de Haan and Piterbarg (2001), tail dependence is described via the measure $\Lambda$ arising as the vague limit in $[0,\infty]^2 \setminus \{(\infty,\infty)\}$ of

$$s \Pr[(s\{1 - F_1(X_1)\}, s\{1 - F_2(X_2)\}) \in \cdot] \xrightarrow{v} \Lambda(\cdot), \qquad s \to \infty. \tag{2.14}$$

Let $P$ be the probability measure on $[0,1]^2$ induced by the random vector $(U_1, U_2) := (1 - F_1(X_1), 1 - F_2(X_2))$. Then (2.14) can be written as

$$t^{-1} P(t\cdot) \xrightarrow{v} \Lambda(\cdot), \qquad t \downarrow 0. \tag{2.15}$$

Comparing (2.14) with (2.2), we find that $\mu$ and $\Lambda$ are connected through a simple change-of-variables formula: for Borel sets $B \subset [0,\infty]^2 \setminus \{(\infty,\infty)\}$,

$$\Lambda(B) = \mu(\{(z_1, z_2) \in \mathbb{E} : (1/z_1, 1/z_2) \in B\}). \tag{2.16}$$

From (2.14) or also from (2.3) and (2.4), it follows that

$$\begin{aligned} \Lambda(c\cdot) &= c\Lambda(\cdot), \qquad 0 < c < \infty, \\ \Lambda([0,u] \times [0,\infty]) &= \Lambda([0,\infty] \times [0,u]) = u, \qquad 0 \le u < \infty. \end{aligned} \tag{2.17}$$

The equality above with $u = 0$ shows that $\Lambda$ does not put any mass on the coordinate axes. Combining (2.6) and (2.16), we find

$$\Phi(\cdot) = \Lambda(\{(u_1, u_2) \in (0,\infty]^2 : \|(u_1^{-1}, u_2^{-1})\| \ge 1, \arctan(u_2/u_1) \in \cdot\}). \tag{2.18}$$

In particular, for $u \in [0,\infty)$,

$$\Lambda(\{\infty\} \times (0,u]) = u\Phi(\{0\}),$$
$$\Lambda((0,u] \times \{\infty\}) = u\Phi(\{\pi/2\}).$$

The spectral measure corresponding to the $L_p$ norm,

$$\|(z_1, z_2)\|_p = \begin{cases} (|z_1|^p + |z_2|^p)^{1/p}, & \text{if } p \in [1,\infty), \\ |z_1| \vee |z_2|, & \text{if } p = \infty, \end{cases}$$

will be denoted by $\Phi_p$. Write

$$y_p(x) = \begin{cases} \infty, & \text{if } x \in [0,1), \\ \left(1 + \dfrac{1}{x^p - 1}\right)^{1/p}, & \text{if } x \in [1,\infty] \text{ and } p \in [1,\infty), \\ 1, & \text{if } x \in [1,\infty] \text{ and } p = \infty. \end{cases} \tag{2.19}$$

Note that for $x \geq 1$, $y_p(x)$ is the (smallest) value of $y \in [1,\infty]$ that solves the equation $\|(x^{-1}, y^{-1})\|_p = 1$. Now by (2.18),

$$\Phi_p([0,\theta]) = \Lambda(C_{p,\theta}), \qquad \theta \in [0,\pi/2],$$

where

$$C_{p,\theta} = \begin{cases} ([0,\infty] \times \{0\}) \cup (\{\infty\} \times [0,1]), \\ \qquad \text{if } \theta = 0, \\ \{(x,y) : 0 \leq x \leq \infty, 0 \leq y \leq (x \tan\theta) \wedge y_p(x)\}, \\ \qquad \text{if } 0 < \theta < \pi/2, \\ \{(x,y) : 0 \leq x \leq \infty, 0 \leq y \leq y_p(x)\}, \\ \qquad \text{if } \theta = \pi/2. \end{cases} \tag{2.20}$$

Further, note that $x \tan\theta < y_p(x)$ if and only if $x < x_p(\theta)$, where for $\theta \in [0,\pi/2]$,

$$x_p(\theta) = \|(1, \cot\theta)\|_p = \begin{cases} (1 + \cot^p \theta)^{1/p}, & \text{if } p \in [1,\infty), \\ 1 \vee \cot\theta, & \text{if } p = \infty. \end{cases} \tag{2.21}$$

The relation between $y_p(x)$, $x_p(\theta)$ and $C_{p,\theta}$ is depicted in Figure 1.

REMARK 2.1. In this paper we shall make no assumptions on the marginal distribution functions $F_1$ and $F_2$ whatsoever except for continuity. However, if in addition to (2.2) the marginal distribution functions $F_1$ and $F_2$ are in the max-domains of attraction of extreme-value distribution functions $G_1$ and $G_2$, then $F$ is actually in the max-domain of attraction of a bivariate extreme-value distribution function $G$ with marginals $G_1$ and $G_2$ and spectral measure $\Phi$. More precisely, if (2.2) holds and if there exist real sequences $a_n > 0$, $b_n$, $c_n > 0$ and $d_n$ such that $F_1^n(a_n x + b_n) \xrightarrow{d} G_1(x)$ and $F_2^n(c_n y + d_n) \xrightarrow{d} G_2(y)$ for all $x, y \in \mathbb{R}$ and as $n \to \infty$, then

$$F^n(a_n x + b_n, c_n y + d_n) \to G(x,y) = \exp[-l\{-\log G_1(x), -\log G_2(y)\}]$$

for all $x, y \in \mathbb{R}$. The function $l$ on the right-hand side is called the *stable tail dependence function* [Huang (1992) and Drees and Huang (1998)] and can be expressed in terms of either $\Lambda$, $\mu$ or $\Phi$ through

$$\begin{aligned} l(x_1, x_2) &= \Lambda(\{(u_1, u_2) \in [0,\infty]^2 : u_1 \le x_1 \text{ or } u_2 \le x_2\}) \\ &= \mu(\{(z_1, z_2) \in [0,\infty]^2 : z_1 \ge x_1^{-1} \text{ or } z_2 \ge x_2^{-1}\}) \\ &= \int_{[0,\pi/2]} \frac{\max(x_1 \sin\theta, x_2 \cos\theta)}{\|(\sin\theta, \cos\theta)\|} \Phi(d\theta), \end{aligned}$$

where we used (2.8) for the final step. Standardizing the marginals of $G$ to the unit-Fréchet distribution yields the extreme-value distribution function

$$G_*(x_1, x_2) = \exp\{-l(1/x_1, 1/x_2)\}, \qquad x_1, x_2 > 0. \tag{2.22}$$

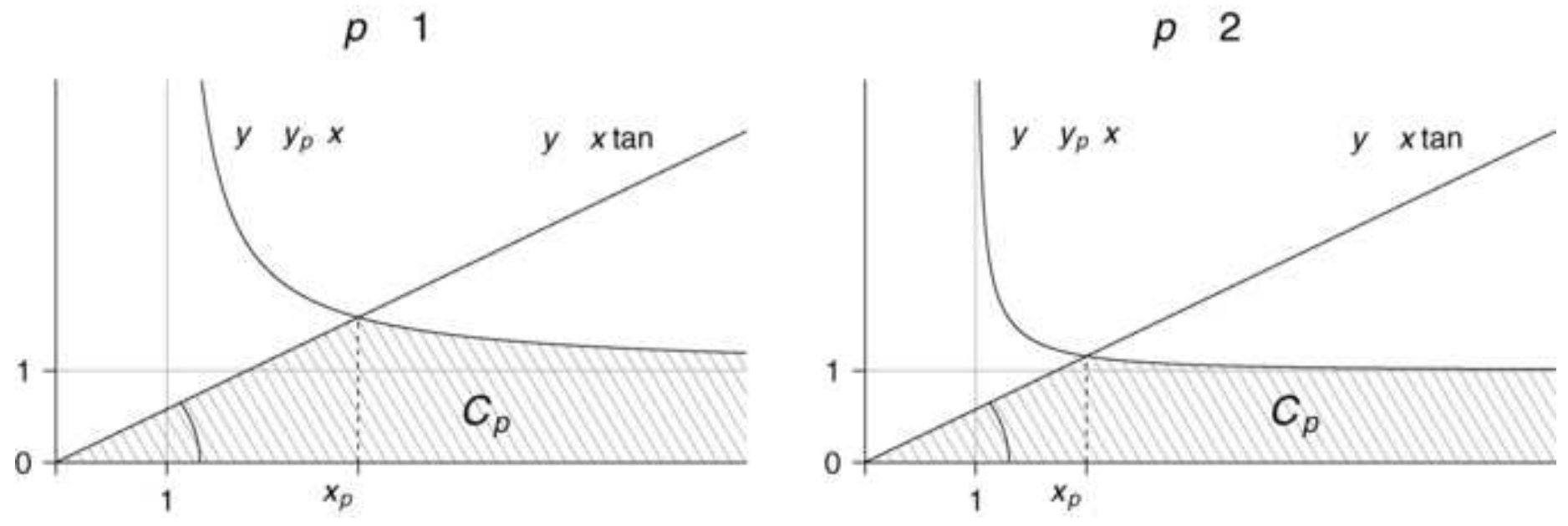

FIG. 1. *The region $C_{p,\theta}$ in (2.20) for $p = 1$ (left) or $p = 2$ (right) and $0 < \theta < \pi/2$.*

In the notation of Coles and Tawn (1991), we have $l(1/x_1, 1/x_2) = V(x_1, x_2)$. The *Pickands dependence function* $A : [0,1] \to [1/2, 1]$ is defined by $A(v) = l(1-v, v)$ for $v \in [0,1]$ [Pickands (1981)]; we have

$$G_*(x_1, x_2) = \exp\left\{-\left(\frac{1}{x_1} + \frac{1}{x_2}\right) A\left(\frac{x_1}{x_1 + x_2}\right)\right\}.$$

The spectral measure $H$ in Coles and Tawn (1991) and Joe, Smith and Weissman (1992) is connected to our $\Phi_1$ through the relation

$$H([0,w]) = \Phi_1(\arctan\{w/(1-w)\}), \qquad w \in [0,1]. \tag{2.23}$$

The Pickands dependence function $A$ can be recovered from $H$ by the formulas

$$\begin{aligned} A(v) &= \int_{[0,1]} \max\{w(1-v), (1-w)v\} H(dw) \\ &= 1 - v + \int_0^v H([0,w])\, dw. \end{aligned} \tag{2.24}$$

Example 2.2 (Asymmetric logistic model). The asymmetric logistic model with parameters $r \ge 1$, $0 \le \psi_1, \psi_2 \le 1$ is defined by its stable tail dependence function

$$l(x_1, x_2) = (1-\psi_1)x_1 + (1-\psi_2)x_2 + \{(\psi_1 x_1)^r + (\psi_2 x_2)^r\}^{1/r}$$

for $(x_1, x_1) \in [0,\infty)^2$ [Tawn (1988)]. The special case $\psi_1 = \psi_2 = 1$ yields the logistic model, originally due to Gumbel (1960). Tail independence arises for $r=1$ or $\psi_1\psi_2=0$, whereas complete tail dependence arises for $\psi_1=\psi_2=1$ and $r\to\infty$. It can be shown that for $1<r<\infty$, the spectral measure $\Phi_p$ has point masses $\Phi_p(\{0\})=1-\psi_2$ and $\Phi_p(\{\pi/2\})=1-\psi_1$, while on the interior $(0,\pi/2)$ the measure $\Phi_p$ is absolutely continuous with density

$$\frac{d\Phi_p}{d\theta}(\theta)=(r-1)(\psi_1\psi_2)^r(\sin^p\theta+\cos^p\theta)^{1/p}(\sin\theta+\cos\theta)^{r-2}$$
$$\times(\sin\theta\cos\theta)^{r-2}\{(\psi_1\cos\theta)^r+(\psi_2\sin\theta)^r\}^{1/r-2}.$$

**3. Empirical spectral measures.** Let $(X_{i1},X_{i2})$, $i=1,\dots,n$, be independent bivariate random vectors from a common distribution function $F$ satisfying (2.2). Our aim is to estimate the spectral measure $\Phi_p$ corresponding to the $L_p$ norm for arbitrary $p\in[1,\infty]$. For convenience, write $\Phi_p(\theta)=\Phi_p([0,\theta])$ for $\theta\in[0,\pi/2]$.

Consider the left-continuous marginal empirical distribution functions

$$\hat F_j(x_j)=\frac{1}{n}\sum_{i=1}^n\mathbf{1}(X_{ij}<x_j),\qquad x_j\in\mathbb{R},\ j=1,2. \tag{3.1}$$

Define

$$\hat U_{ij}=1-\hat F_j(X_{ij})=\frac{n+1-R_{ij}}{n},\qquad i=1,\dots,n;\ j=1,2; \tag{3.2}$$

here $R_{ij}=\sum_{l=1}^n\mathbf{1}(X_{lj}\le X_{ij})$ is the rank of $X_{ij}$ among $X_{1j},\dots,X_{nj}$. Let $\hat P_n$ be the empirical measure of $(\hat U_{i1},\hat U_{i2})$, $i=1,\dots,n$, that is,

$$\hat P_n(\cdot)=\frac{1}{n}\sum_{i=1}^n\mathbf{1}\{(\hat U_{i1},\hat U_{i2})\in\cdot\}.$$

Observe that the transformed data $(\hat U_{i1},\hat U_{i2})$, $i=1,\dots,n$, are no longer independent. This dependence will contribute to the limiting distribution of the estimators to be considered.

Let $k=k_n\in(0,n]$ be an intermediate sequence, that is, $k\to\infty$ and $k/n\to 0$ as $n\to\infty$. We find our estimator $\hat\Phi_p$ by using (2.15) and (2.18) with $t=k/n$ and $P$ replaced by $\hat P_n$. In terms of distribution functions, this becomes

$$\begin{aligned}\hat\Phi_p(\theta)&=\frac{n}{k}\hat P_n\left(\frac{k}{n}C_{p,\theta}\right)\\&=\frac{1}{k}\sum_{i=1}^n\mathbf{1}\{(n+1-R_{i1})^{-p}+(n+1-R_{i2})^{-p}\ge k^{-p},\\&\qquad\qquad n+1-R_{i2}\le(n+1-R_{i1})\tan\theta\}\end{aligned}$$

for $\theta\in[0,\pi/2]$ and with $C_{p,\theta}$ as in (2.20).

In Einmahl, de Haan and Piterbarg (2001), the limiting behavior of $\hat{\Phi}_p$ has been derived in case $p=\infty$. We now present a generalization to all $L_p$ norms for $p \in [1,\infty]$. More precisely, we will study the asymptotic behavior of the process

$$\sqrt{k}\{\hat{\Phi}_p(\theta)-\Phi_p(\theta)\}, \qquad \theta \in [0,\pi/2].$$

We will assume that

$$\Lambda = \Lambda_c + \Lambda_d, \tag{3.3}$$

where $\Lambda_c$ is absolutely continuous with a density $\lambda$, which is continuous on $[0,\infty)^2 \setminus \{(0,0)\}$, and with $\Lambda_d$ such that $\Lambda_d([0,\infty)^2)=0$, $\Lambda_d(\{\infty\}\times[0,u]) = u\Phi_p(\{0\})$ and $\Lambda_d([0,u]\times\{\infty\}) = u\Phi_p(\{\pi/2\})$ for $u \in [0,\infty)$. In contrast to Einmahl, de Haan and Piterbarg (2001), $\Phi_p$ is allowed to have atoms at 0 and $\pi/2$; in particular, tail independence is allowed. Also, the restriction of $\Phi_p$ to $(0,\pi/2)$ is absolutely continuous with a continuous density. This excludes complete tail dependence, that is, $\Phi_p$ being degenerate at $\pi/4$, in which case $\Lambda$ is concentrated on the diagonal. The homogeneity of $\Lambda$ in (2.17) implies that $\lambda(cu_1,cu_2) = c^{-1}\lambda(u_1,u_2)$ for all $c>0$ and $(u_1,u_2) \in [0,\infty)^2 \setminus \{(0,0)\}$.

Let $P_n$ be the empirical measure of $(U_{i1},U_{i2}) = (1-F_1(X_{i1}), 1-F_2(X_{i2}))$, $i=1,\dots,n$, and let $\Gamma_{jn}(u) = n^{-1}\sum_{i=1}^n \mathbf{1}(U_{ij}\le u)$, $u\in[0,1]$ and $j\in\{1,2\}$, be the corresponding marginal empirical distribution functions; for $u \in (1,\infty]$, we set $\Gamma_{jn}(u)=u$. Furthermore, for $\theta\in[0,\pi/2]$, define the set

$$\hat{C}_{p,\theta} = \frac{n}{k}\left\{(u_1,u_2)\in[0,\infty]^2\setminus\{(\infty,\infty)\} : (\Gamma_{1n}(u_1),\Gamma_{2n}(u_2)) \in \frac{k}{n}C_{p,\theta}\right\}.$$

From the identity $\Gamma_{jn}(u) = 1-\hat{F}_j(F_j^{-1}(1-u))$ for $u\in(0,1)$ it follows that

$$\hat{P}_n\left(\frac{k}{n}C_{p,\theta}\right) = P_n\left(\frac{k}{n}\hat{C}_{p,\theta}\right).$$

This representation yields the following crucial decomposition: for $\theta\in[0,\pi/2]$,

$$\begin{aligned}
\sqrt{k}\{\hat{\Phi}_p(\theta)-\Phi_p(\theta)\} &= \sqrt{k}\left\{\frac{n}{k}P_n\left(\frac{k}{n}\hat{C}_{p,\theta}\right) - \frac{n}{k}P\left(\frac{k}{n}\hat{C}_{p,\theta}\right)\right\}\\
&\quad + \sqrt{k}\left\{\frac{n}{k}P\left(\frac{k}{n}\hat{C}_{p,\theta}\right) - \Lambda(\hat{C}_{p,\theta})\right\}\\
&\quad + \sqrt{k}\{\Lambda(\hat{C}_{p,\theta}) - \Lambda(C_{p,\theta})\}\\
&=: V_{n,p}(\theta) + r_{n,p}(\theta) + Y_{n,p}(\theta).
\end{aligned} \tag{3.4}$$

The first term, $V_{n,p}$, features a local empirical process evaluated in a random set $\hat{C}_{p,\theta}$. The second term, $r_{n,p}$, is a bias term, which will vanish in the limit under our assumptions. The third term, $Y_{n,p}$, is due to the fact that the marginal distributions are unknown and captures the effect of the rank transformation in (3.1) and (3.2).

Next we will define the processes that will arise as the weak limits of the processes $V_{n,p}$ and $Y_{n,p}$ in (3.4). Define $W_\Lambda$ to be a Wiener process indexed by the Borel sets of $[0,\infty]^2 \setminus \{(\infty,\infty)\}$ and with "time" $\Lambda$, that is, a centered Gaussian process with covariance function $E[W_\Lambda(C)W_\Lambda(C')] = \Lambda(C \cap C')$. We can write, in the obvious notation, $W_\Lambda = W_{\Lambda_c} + W_{\Lambda_d}$, where the two processes on the right are independent. Note that

$$(W_\Lambda(C_{p,\theta}))_{\theta\in[0,\pi/2]} \stackrel{d}{=} (W(\Phi_p(\theta)))_{\theta\in[0,\pi/2]}$$

with $W$ a standard Wiener process on $[0,\infty)$. Define $W_1(x) = W_\Lambda([0,x] \times [0,\infty])$ and $W_2(y) = W_\Lambda([0,\infty] \times [0,y])$ for $x, y \in [0,\infty)$. Note that $W_1$ and $W_2$ are standard Wiener processes as well. For $p \in [1,\infty)$, define the process $Z_{c,p}$ on $[0,\pi/2]$ by

$$\begin{aligned} Z_{c,p}(\theta) &= \mathbf{1}(\theta < \pi/2) \int_0^{x_p(\theta)} \lambda(x, x\tan\theta)\{W_1(x)\tan\theta - W_2(x\tan\theta)\}\,dx \\ &\quad + \begin{cases} \displaystyle\int_{x_p(\theta)}^\infty \lambda(x, y_p(x))\{W_1(x)y_p'(x) - W_2(y_p(x))\}\,dx, \\ \qquad \text{if } p < \infty, \\ \displaystyle -W_1(1)\int_1^{1\vee\tan\theta} \lambda(1,y)\,dy - W_2(1)\int_{1\vee\cot\theta}^\infty \lambda(x,1)\,dx, \\ \qquad \text{if } p = \infty \end{cases} \end{aligned}$$

with $y_p'$ the derivative of $y_p$. Define $Z_d$ by

$$Z_d(\theta) = -\Phi_p(\{0\})W_2(1), \qquad \theta \in [0,\pi/2],$$

and write $Z_p = Z_{c,p} + Z_d$. It is our aim to show that

$$(V_{n,p}, r_{n,p}, Y_{n,p}) \xrightarrow{d} (W_\Lambda(C_{p,\cdot}), 0, Z_p), \qquad n \to \infty.$$

This convergence and the decomposition in (3.4) then will yield the asymptotic behavior of $\sqrt{k}(\hat\Phi_p - \Phi)$.

Assume that $P$ is absolutely continuous with density $p$. Then the measure $t^{-1}P(t\cdot)$, for $t > 0$, is absolutely continuous as well with density $tp(tu_1, tu_2)$. For $1 \le T < \infty$ and $t > 0$, define

$$\mathcal{D}_T(t) := \iint_{\mathcal{L}_T} |tp(tu_1,tu_2) - \lambda(u_1,u_2)|\,du_1\,du_2, \tag{3.5}$$

where $\mathcal{L}_T = \{(u_1,u_2) : 0 \le u_1 \wedge u_2 \le 1, u_1 \vee u_2 \le T\}$.

THEOREM 3.1. *Assume the framework of Section 2 and suppose $\Lambda$ is as in (3.3). Then, if $\mathcal{D}_{1/t}(t) \to 0$ as $t \downarrow 0$ and if the intermediate sequence $k$ is such that*

$$\sqrt{k}\mathcal{D}_{n/k}(k/n) \to 0, \qquad n \to \infty, \tag{3.6}$$

*then in $D[0, \pi/2]$ and as $n \to \infty$,*

$$\sqrt{k}(\hat{\Phi}_p - \Phi_p) \xrightarrow{d} W_\Lambda(C_{p,\cdot}) + Z_p =: \alpha_p. \tag{3.7}$$

The condition $\lim_{t\downarrow 0} \mathcal{D}_{1/t}(t) = 0$ in Theorem 3.1 implies $\Phi_p(\{0, \pi/2\}) = 0$ and thus $\Lambda_d = 0$. Indeed, in case $\Lambda_d \neq 0$, the convergence in (3.7) cannot hold: when, for example, $\Phi_p(\{0\}) > 0$, we have, since $\hat{\Phi}_p(0) = 0$, $\sqrt{k}\{\hat{\Phi}_p(0) - \Phi_p(0)\} \to -\infty$. In contrast, the following result does allow $\Phi_p$ to have atoms at 0 or $\pi/2$. Recall $\mathcal{D}_T(t)$ in (3.5) and $\alpha_p$ in (3.7).

THEOREM 3.2. *Let $\eta \in (0, \pi/4)$. Assume the framework of Section 2 and suppose $\Lambda$ is as in (3.3). Then, if $\mathcal{D}_1(t) \to 0$ as $t \downarrow 0$ and if the intermediate sequence $k$ is such that*

$$\sqrt{k} \inf_{T>0} \{\mathcal{D}_T(k/n) + 1/T\} \to 0, \qquad n \to \infty, \tag{3.8}$$

*then in $D[\eta, \pi/2 - \eta]$ and as $n \to \infty$,*

$$\sqrt{k}(\hat{\Phi}_p - \Phi_p) \xrightarrow{d} \alpha_p. \tag{3.9}$$

In case of tail independence, that is, $\Phi_p(\{0\}) = \Phi_p(\{\pi/2\}) = 1$ and $\lambda = 0$, we have $\alpha_p = 0$.

Under a stronger condition on the sequence $k$, the convergence of the process $\sqrt{k}(\hat{\Phi}_p - \Phi_p)$ holds on the whole interval $[0, \pi/2]$, provided that we flatten the process on intervals $[0, \eta_n]$ and $[\pi/2 - \eta_n, \pi/2]$, with $\eta_n \in (0, \pi/4)$ tending to zero sufficiently slowly. Define the transformation $\tau_n : [0, \pi/2] \to [0, \pi/2]$ by

$$\tau_n(\theta) = \begin{cases} \eta_n, & \text{if } 0 \le \theta < \eta_n, \\ \theta, & \text{if } \eta_n \le \theta \le \pi/2 - \eta_n, \\ \pi/2 - \eta_n, & \text{if } \eta_n < \theta < \pi/2, \\ \pi/2, & \text{if } \theta = \pi/2. \end{cases} \tag{3.10}$$

THEOREM 3.3. *Let $k$ be an intermediate sequence and let $\eta_n = (k/n)^a$ for some fixed $a \in (0, 1)$. Assume the framework of Section 2 and suppose $\Lambda$ is as in (3.3). If*

$$\sqrt{k} \inf_{T \ge 2/\eta_n} \{\mathcal{D}_T(k/n) + 1/T\} \to 0, \qquad n \to \infty, \tag{3.11}$$

*then in* $D[0,\pi/2]$ *and as* $n\to\infty$,

$$\sqrt{k}(\hat{\Phi}_p-\Phi_p)\circ\tau_n \xrightarrow{d} \alpha_p. \tag{3.12}$$

Theorems 3.1 and 3.3 will be instrumental when establishing our main results in the next section.

**4. Enforcing the moment constraints.** Fix $p\in[1,\infty]$ and let $\mathcal{Q}_p$ be the class of probability measures $Q_p$ on $[0,\pi/2]$ such that

$$\int_{[0,\pi/2]} f(\theta)Q_p(d\theta)=0, \tag{4.1}$$

where

$$f(\theta)=f_p(\theta)=\frac{\sin\theta-\cos\theta}{\|(\sin\theta,\cos\theta)\|_p}, \qquad \theta\in[0,\pi/2]. \tag{4.2}$$

If $Q_p$ is the spectral probability measure of some exponent measure $\mu$ with respect to the $L_p$ norm, then $Q_p\in\mathcal{Q}_p$ by (2.12). Conversely, if $Q_p\in\mathcal{Q}_p$, then we can define an exponent measure $\mu$ through (2.8) and (2.13) which has $Q_p$ as its spectral probability measure with respect to the $L_p$ norm. As before, denote distribution functions of measures under consideration by $Q_p(\theta)=Q_p([0,\theta])$, etc.

In view of (2.10), we define the empirical spectral probability measure $\hat{Q}_p$ by

$$\hat{Q}_p(\cdot)=\frac{\hat{\Phi}(\cdot)}{\hat{\Phi}(\pi/2)}=\frac{1}{N_n}\sum_{i\in I_n}\mathbf{1}(\Theta_{in}\in\cdot), \tag{4.3}$$

where $N_n=|I_n|$ and

$$\Theta_{in}=\arctan(\hat{U}_{i2}/\hat{U}_{i1}), \qquad i=1,\ldots,n;$$
$$I_n=\{i=1,\ldots,n:\|(\hat{U}_{i1}^{-1},\hat{U}_{i2}^{-1})\|_p\geq n/k\}.$$

Typically,

$$\int f\,d\hat{Q}_p=\frac{1}{N_n}\sum_{i\in I_n}f(\Theta_{in})$$

is different from zero, in which case $\hat{Q}_p$ does not belong to $\mathcal{Q}_p$, that is, $\hat{Q}_p$ is itself not a spectral probability measure.

Therefore, we propose to modify $\hat{Q}_p$ such that the moment constraint (4.1) is fulfilled and the new estimator does belong to $\mathcal{Q}_p$: define

$$\tilde{Q}_p(\cdot):=\sum_{i\in I_n}\tilde{p}_{in}\mathbf{1}(\Theta_{in}\in\cdot),$$

where the weight vector $(\tilde{p}_{in} : i \in I_n)$ solves the following optimization problem:

$$
\begin{aligned}
&\text{maximize} && \prod_i p_{in}, \\
&\text{constraints} && p_{in} \geq 0 \qquad \text{for all } i \in I_n, \\
&&& \sum_i p_{in} = 1, \\
&&& \sum_i p_{in} f(\Theta_{in}) = 0.
\end{aligned}
\tag{4.4}
$$

The thus obtained estimator $\tilde{Q}_p$ can be viewed as a maximum empirical likelihood estimator (MELE) based on the sample $\{\Theta_{in} : i \in I_n\}$; see the monograph Owen (2001). Actually, the optimization problem in (4.4) can be readily solved by the method of Lagrange multipliers [see, e.g., Owen (2001), page 22]: let $\tilde{\mu}_n$ be the solution in $(-1, 1)$ to the nonlinear equation

$$
\sum_{i \in I_n} \frac{f(\Theta_{in})}{1 + \tilde{\mu}_n f(\Theta_{in})} = 0 \tag{4.5}
$$

and define

$$
\tilde{p}_{in} = \frac{1}{N_n} \frac{1}{1 + \tilde{\mu}_n f(\Theta_{in})}, \qquad i \in I_n, \tag{4.6}
$$

then the vector $(\tilde{p}_{in} : i \in I_n)$ is the solution to (4.4). Observe that the original estimator $\hat{Q}_p$ corresponds to $\tilde{\mu}_n = 0$ and is the solution to (4.4) without the final constraint $\sum_i p_{in} f(\Theta_{in}) = 0$.

Since $\tilde{Q}_p \in \mathcal{Q}_p$, we can exploit the transformation formulas in Section 2 to define estimators of the spectral measure $\Phi_p$: as in (2.13),

$$
\tilde{\Phi}_p(\cdot) := \tilde{Q}_p(\cdot)/m_p(\tilde{Q}_p),
$$

where for a bounded, measurable function $h : [0, \pi/2] \to \mathbb{R}$,

$$
m_p(h) := -\int_0^{\pi/2} h(\theta)\, d\frac{\cos\theta}{\|(\sin\theta, \cos\theta)\|_p};
$$

cf. (2.12). Further, for $\theta \in [0, \pi/2]$, define $I(\theta) = \int_{[0,\theta]} f(\vartheta)\, dQ_p(\vartheta)$ and

$$
\begin{aligned}
\beta_p(\theta) &= \frac{\Phi_p(\pi/2)\alpha_p(\theta) - \alpha_p(\pi/2)\Phi_p(\theta)}{\Phi_p^2(\pi/2)}, \\
\gamma_p(\theta) &= \beta_p(\theta) + \frac{\int_{[0,\pi/2]} \beta_p\, df}{\int_{[0,\pi/2]} f^2\, dQ_p} I(\theta), \\
\delta_p(\theta) &= \frac{m_p(Q_p)\gamma_p(\theta) - m_p(\gamma_p) Q_p(\theta)}{m_p^2(Q_p)}.
\end{aligned}
\tag{4.7}
$$

Note that under the assumptions of Theorem 4.1 below, $Q_p(\{\pi/4\}) < 1$ and thus $\int f^2 \, dQ_p > 0$, so that $\gamma_p(\theta)$ is well defined.

The next two theorems, providing asymptotic normality of $\tilde{\Phi}_p$, are the main results of this paper.

THEOREM 4.1. *Let the assumptions of Theorem 3.1 be fulfilled. Then with probability tending to one, (4.5) admits a unique solution $\tilde{\mu}_n$ and hence in this case the vector $(\tilde{p}_{in} : i \in I_n)$ in (4.6) is the unique solution to (4.4). Also, in $D[0, \pi/2]$ and as $n \to \infty$,*

$$\sqrt{k}(\tilde{Q}_p - Q_p) \xrightarrow{d} \gamma_p, \tag{4.8}$$

$$\sqrt{k}(\tilde{\Phi}_p - \Phi_p) \xrightarrow{d} \delta_p. \tag{4.9}$$

Since Theorem 4.1 is based on Theorem 3.1, the spectral measure cannot have atoms at 0 or $\pi/2$. The following result, based on Theorem 3.3, does allow for such atoms.

THEOREM 4.2. *Fix $\eta \in (0, \pi/4)$ and let $\eta_n = (k/n)^a$ for some $0 < a < 1$. Assume the framework of Section 2 and suppose $\Lambda$ is as in (3.3). If*

$$\sqrt{k}\mathcal{D}_{2/\eta_n}(k/n) + \sqrt{k}\eta_n \to 0, \qquad n \to \infty, \tag{4.10}$$

*then in $D[\eta, \pi/2 - \eta]$ and as $n \to \infty$, the convergence in (4.8) and (4.9) holds.*

REMARK 4.3. From (4.7), it is straightforward to express the limit process $\delta_p$ in terms of the process $\alpha_p$ and thus of $W_\Lambda$. However, because of the presence of the process $Z_p$, no major simplification occurs. As a consequence, we were not able to show that $\tilde{\Phi}_p$ is asymptotically more efficient than $\hat{\Phi}_p$. However, the simulation study in Section 5 does indicate that enforcing the moment constraints leads to a sizeable improvement of the estimator's performance.

REMARK 4.4. Replacing $\Phi_1$ by $\tilde{\Phi}_1$ in (2.23)–(2.24) yields an estimator $\tilde{A}$ of the Pickands dependence function $A$ that is itself a genuine Pickands dependence function. The weak limit of the process $\sqrt{k}(\tilde{A} - A)$ in the function space $C[0,1]$ can be easily derived from the one of $\sqrt{k}(\tilde{\Phi}_1 - \Phi_1)$. Nonparametric estimation of a Pickands dependence function in the domain-of-attraction context was also studied in Capéraà and Fougères (2000) and Abdous and Ghoudi (2005).

## 5. Examples, simulations and real data analyses.

### 5.1. *Examples.*

EXAMPLE 5.1 (Cauchy). Consider the bivariate Cauchy distribution on $(0,\infty)^2$ with density $(2/\pi)(1+x^2+y^2)^{-3/2}$ for $x,y>0$. It follows that

$$\Lambda([0,x]\times[0,y]) = x+y-(x^2+y^2)^{1/2}, \qquad 0\le x,y<\infty$$

and

$$\Phi_p(\theta) = \int_0^\theta \|(\sin\vartheta,\cos\vartheta)\|_p\,d\vartheta$$

for $\theta\in[0,\pi/2]$. It can be shown that $\mathcal{D}_{1/t}(t)=O(t)$ as $t\downarrow 0$. Therefore, Theorems 3.1 and 4.1 hold when $k=o(n^{2/3})$ as $n\to\infty$. In Figure 2, the function $\Phi_p$ is plotted for $p\in\{1,2,3,\infty\}$ together with trajectories of the empirical spectral measure and the MELE for a single sample of size $n=1000$ and at $k=30$. Observe that both estimators share the same set of atoms but with possibly different weights.

The bivariate Cauchy distribution on the full plane $\mathbb{R}^2$ has density function $(2\pi)^{-1}(1+x^2+y^2)^{-3/2}$ for $x,y\in\mathbb{R}$ and spectral measure

$$\Phi_p(\theta) = \frac{1}{2}\left(1+\int_0^\theta \|(\sin\vartheta,\cos\vartheta)\|_p\,d\vartheta + \mathbf{1}(\theta=\pi/2)\right).$$

In particular, $\Phi_p(\{0\})=\Phi_p(\{\pi/2\})=1/2$. For every $0<a<1$ and $\eta_n=(k/n)^a$, we find $\mathcal{D}_{2/\eta_n}(k/n)=O((k/n)^{2-a})$ as $n\to\infty$. Hence the conclusions of Theorems 3.3 and 4.2 hold provided $k=o(n^{2a/(2a+1)})$ as $n\to\infty$. In fact, the results of Theorem 4.2 can be shown to hold when $k=o(n^{2/3})$ as $n\to\infty$.

EXAMPLE 5.2 (Mixture). For $r\in[0,1]$, consider the bivariate distribution function

$$F(x,y) = \left(1-\frac{1}{x}\right)\left(1-\frac{1}{y}\right)\left(1+\frac{r}{x+y}\right), \qquad x,y\ge 1;$$

cf. de Haan and Resnick (1977), Example 3. Its density can be written as a mixture of two densities, $(1-r)f_1(x,y)+rf_2(x,y)$, where

$$f_1(x,y)=\frac{1}{x^2y^2}, \qquad f_2(x,y)=\frac{2}{(x+y)^3}\left(1+\frac{x^2+3xy+y^2}{x^2y^2}\right), \qquad x,y\ge 1.$$

Note that $f_1$ is the density of two independent Pareto(1) random variables. Obviously for $r=0$ we have (tail) independence. The law $P$ of $(1-F_1(X),1-F_2(Y))=(1/X,1/Y)$ is determined by

$$P([0,u]\times[0,v]) = uv\left(1+r\frac{(1-u)(1-v)}{u+v}\right), \qquad 0<u,v\le 1,$$

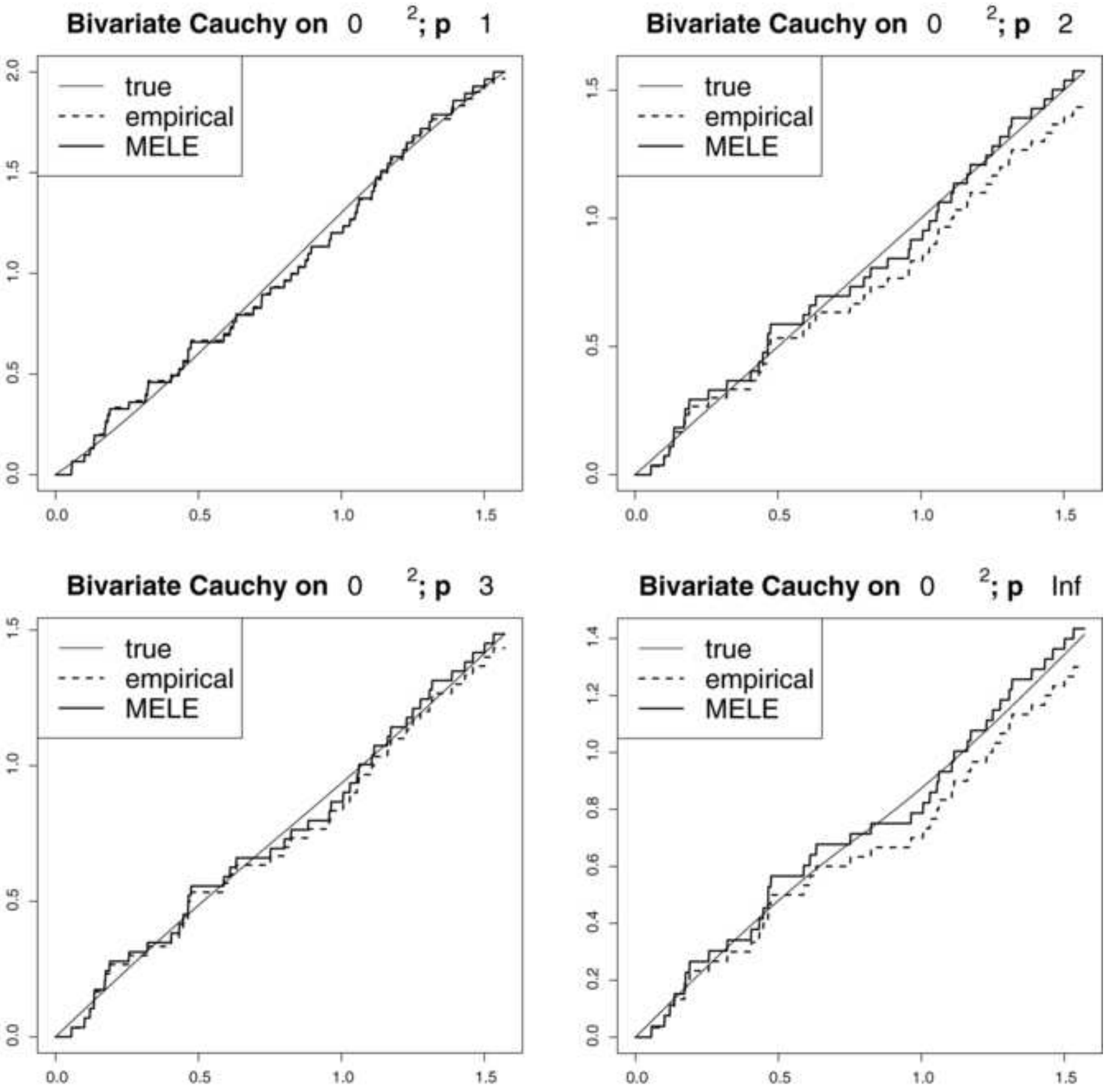

FIG. 2. *Trajectories of the empirical spectral measure (dashed) and the MELE (solid) for a single sample of size $n = 1000$ and at $k = 30$ from the bivariate Cauchy distribution on $(0,\infty)^2$ and for $p \in \{1, 2, 3, \infty\}$.*

hence

$$\Lambda([0,x] \times [0,y]) = r\frac{xy}{x+y}, \qquad 0 < x, y < \infty.$$

For $p \in [1,\infty]$, the corresponding spectral measure $\Phi_p$ satisfies $\Phi_p(\{0\}) = \Phi_p(\{\pi/2\}) = 1 - r$ and

$$\Phi_p(\theta) = 1 - r + 2r\int_0^\theta \frac{\|(\sin\vartheta, \cos\vartheta)\|_p}{(\cos\vartheta + \sin\vartheta)^3}\,d\vartheta + (1-r)\mathbf{1}(\theta = \pi/2)$$

for $\theta \in [0, \pi/2]$. It can be seen that $\mathcal{D}_T(t) = TO(t)$ as $t \downarrow 0$, uniformly in $T > 0$. As a consequence, conditions (3.11) and (4.10) in Theorems 3.3 and 4.2 hold for $a = 1/2$ provided $k = o(n^{1/2})$ as $n \to \infty$. If $r = 1$, the spectral measure $\Phi_p$ has no atoms. Then $\mathcal{D}_{1/t}(t) = O(t)$ as $t \downarrow 0$, so that condition (3.6) in Theorems 3.1 and 4.1 holds provided $k = o(n^{2/3})$ as $n \to \infty$.

5.2. *Simulations.* In order to compare the two spectral measure estimators, a simulation study was conducted involving the following distributions:

- Figure 3: the extreme-value distribution $G_*$ in (2.22) with logistic dependence structure $r = 2$ and $\psi_1 = \psi_2 = 1$ (Example 2.2);
- Figure 4: the bivariate Cauchy distribution on $(0, \infty)^2$ (Example 5.1);
- Figure 5: the mixture distribution with $r = 0.5$ (Example 5.2).

For each distribution, 1000 samples of size $n = 1000$ were drawn. For each sample, the empirical spectral measure and the MELE were computed for various ranges of $k$ and at $p \in \{1, 2, 3, \infty\}$. For each estimate, the Integrated Squared Errors $\int_a^b (\hat{\Phi}_p - \Phi_p)^2$ and $\int_a^b (\tilde{\Phi}_p - \Phi_p)^2$ were computed, where $(a, b) = (0, \pi/2)$ for the logistic model and the Cauchy distribution on $(0, \infty)^2$ whereas $(a, b) = (0.05\pi/2, 0.95\pi/2)$ for the mixture model. Next, the ISEs were averaged out over the 1000 samples, yielding the estimates of the Mean Integrated Squared Errors which are displayed in the figures.

From the plots, it can be seen that in all these cases the minimum MISE of the MELE is smaller than the one of the empirical spectral measure. Overall the MELE outperforms the empirical spectral measure, but in two cases, away from the minimum MISE, the empirical spectral measure performs slightly better than the MELE (Figure 5). Moreover, most of the time, the MISE is decreasing in $p$. Finally, the presence of atoms at the endpoints has an adverse effect on the estimators.

5.3. *Case studies.*

5.3.1. *Loss-ALAE data.* In Frees and Valdez (1998), copula fitting was illustrated on a dataset of 1500 insurance company indemnity claims, displayed in Figure 6. Each claim consists of an indemnity payment (Loss) and an allocated loss adjustment expense (ALAE), that is, an expense specifically attributable to the settlement of the individual claim such as lawyer's fees and investigation expenses.

The problem is which bivariate tail dependence structure to fit to these data. The first question is whether the data are tail independent or not. If yes, then the spectral measure should be the sum of the Dirac measures at 0 and $\pi/2$. In Figure 6, the MELE ($p = 1$, $k = 40$) shows substantial increase away from the endpoints, in clear contrast to the case of independence, for which the spectral measure is flat on $(0, \pi/2)$. Note that a formal test for tail independence is not possible, since under the null hypothesis of tail independence the limit process $\alpha_p$ in Theorem 3.1 (and hence all limit processes in subsequent theorems) is degenerate at zero; in contrast, tests with tail independence in the *alternative* hypothesis can be constructed in the more specialized framework of Ledford and Tawn (1996). By way of illustration, we compared the MELE for the Loss-ALAE data with the MELE for a random sample of the same size from the bivariate lognormal distribution fitted to the data (right-hand panel of Figure 6). The MELE for this sample is compared to the MELE of the Loss-ALAE data in Figure 7 (left-hand panel). Despite the fact that the correlation (on log scale) of both samples is equal to 0.44, the MELE of the lognormal sample is much closer to the spectral measure of independence (which is the true one for this sample).

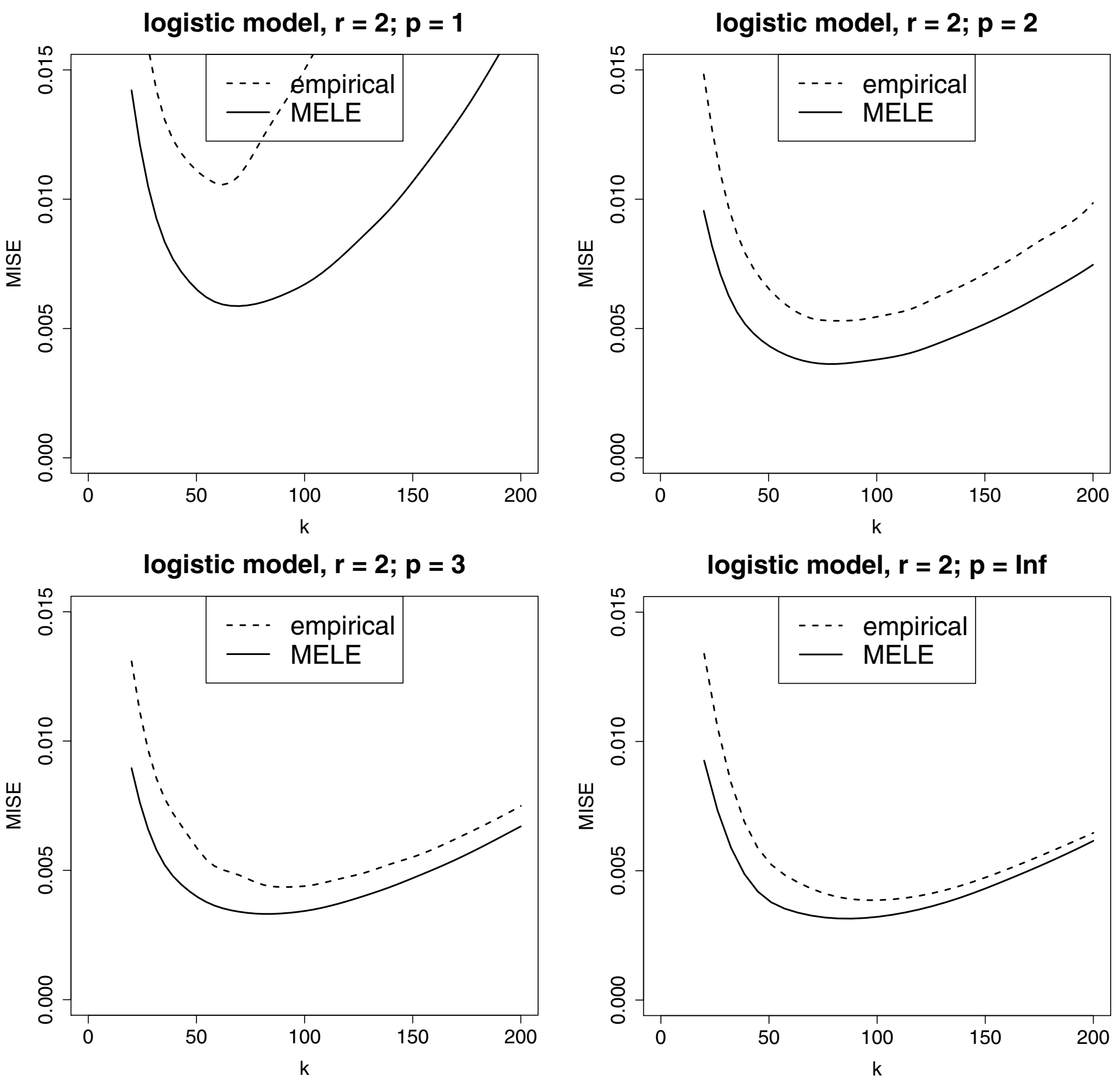

Fig. 3. *MISE of the empirical spectral measure (dashed) and the MELE (solid) for* 1000 *samples of size* $n = 1000$ *from the logistic model with* $r = 2$ *and* $p \in \{1, 2, 3, \infty\}$.

The next question is then which model to fit to the data. In Genest, Ghoudi and Rivest (1998), the logistic model (Example 2.2) is pro-posed, while in Beirlant et al. (2004), Chapter 9, both the logistic ($r = 0.73$) and the asymmetric logistic model ($r = 0.66$, $\psi_1 = 1$, $\psi_2 = 0.89$) are fitted. In the left-hand panel of Figure 7, the spectral measures ($p = 1$) of both parametric models are displayed and compared to the MELE. The asymmetric logistic model follows the MELE much more closely. Interestingly, the spectral measure of the fitted asymmetric logistic model does not have an atom at $\pi/2$ while it has an atom at 0 of size $1 - \psi_2 = 0.11$. Formal goodness-of-fit tests based on the stable tail dependence function $l$ are described in de Haan, Neves and Peng (2008) and Einmahl, Krajina and Segers (2008).

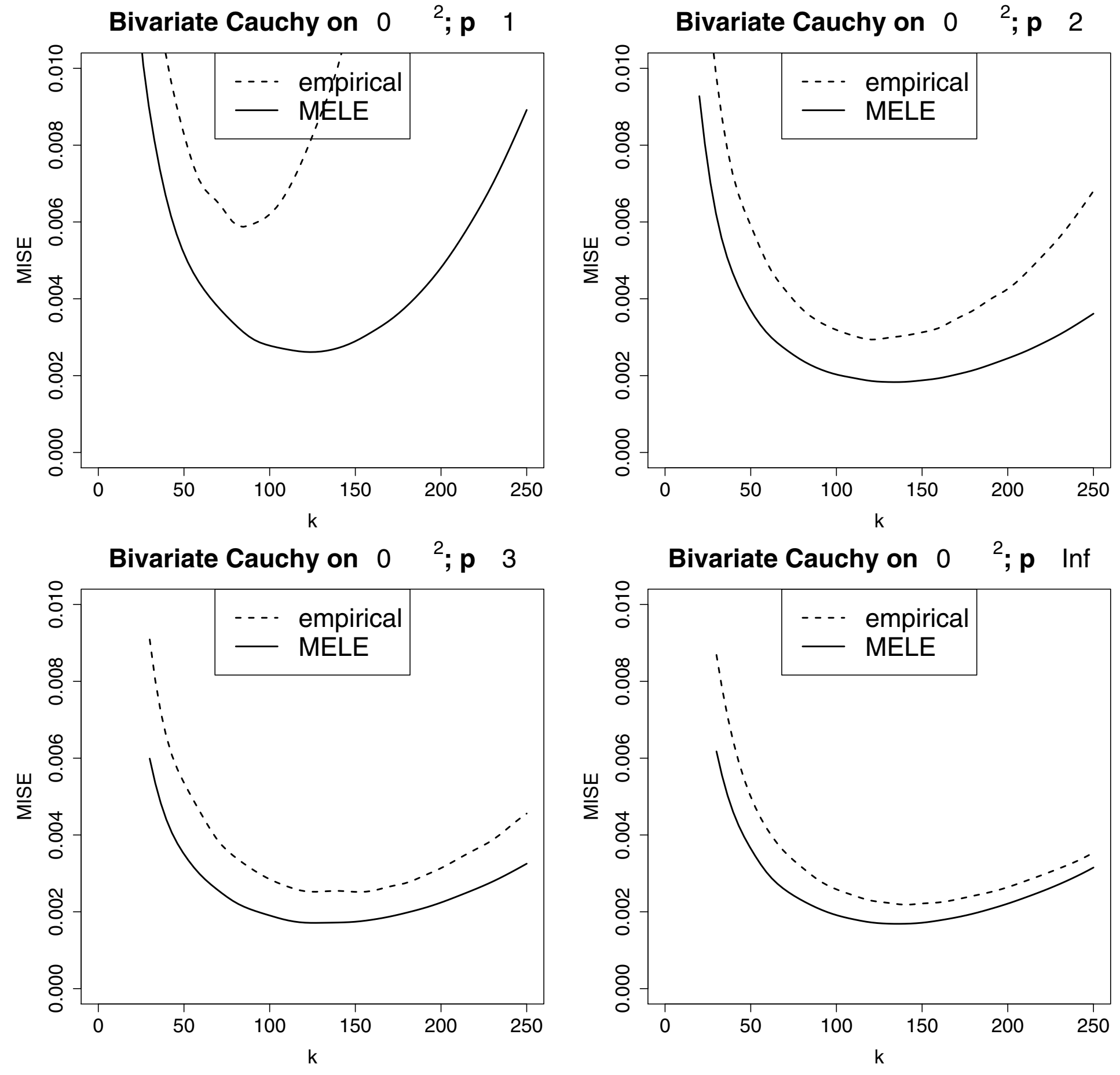

FIG. 4. *MISE of the empirical spectral measure (dashed) and the MELE (solid) for* 1000 *samples of size* $n = 1000$ *from the bivariate Cauchy distribution on the positive quadrant and* $p \in \{1, 2, 3, \infty\}$.

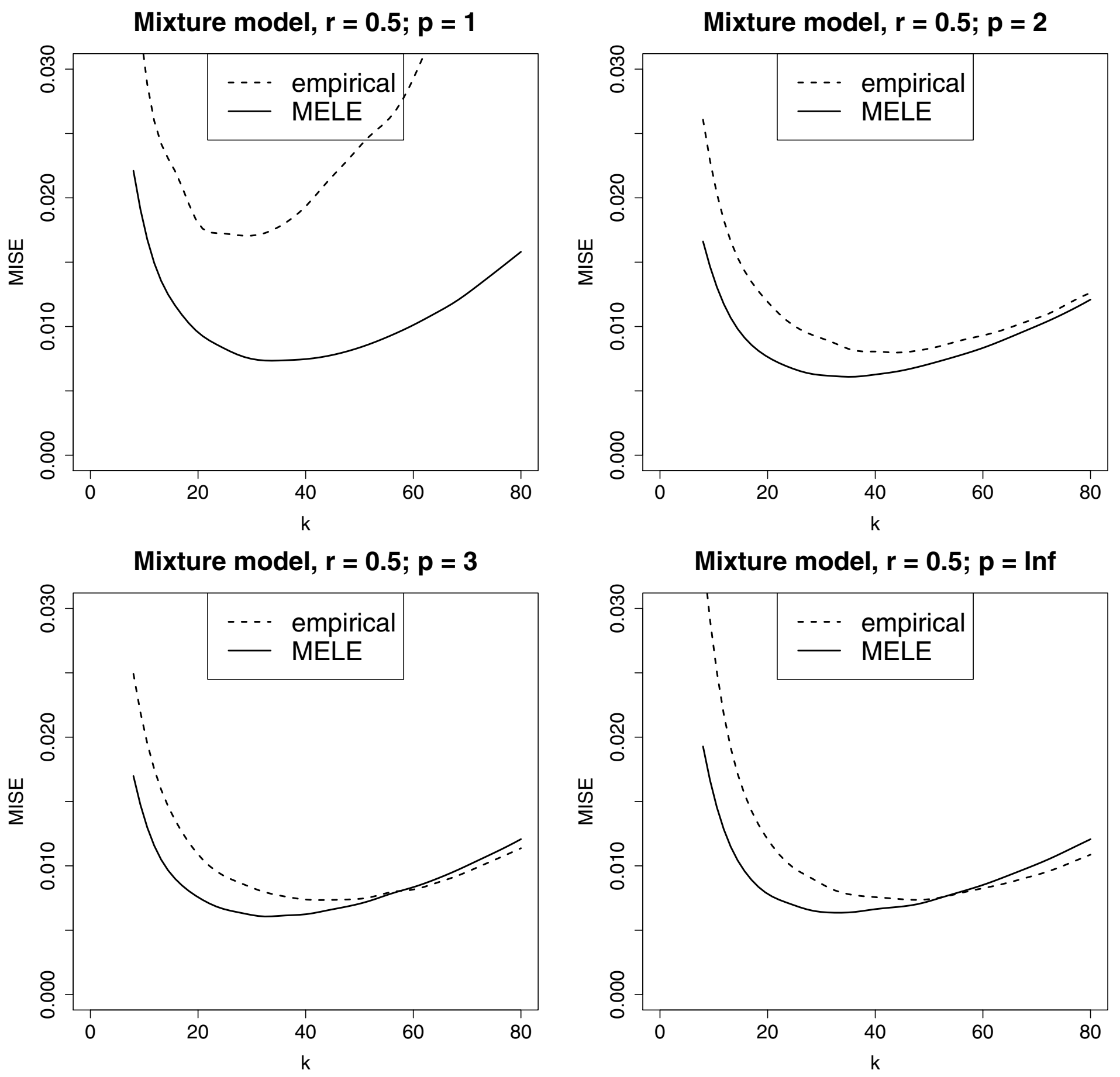

FIG. 5. *MISE of the empirical spectral measure (dashed) and the MELE (solid) for* 1000 *samples of size* $n = 1000$ *from the mixture distribution with* $r = 0.5$ *and* $p \in \{1, 2, 3, \infty\}$.

5.3.2. *NHANES data.* We consider the variables Standing Height (cm) and Weight (kg) for females from the National Health and Nutrition Examination Survey (NHANES) 2005–2006 [National Center for Health Statistics (2007)]. Data of females of age 18–64 years only have been retained; this leads to a sample size of 2237. In Figure 8, the MELE of the spectral measure ($p = 1$) has been computed for $k \in \{25, 50\}$ and compared to the spectral measure for tail independence. Although the data have a nonnegligible positive correlation (0.283), the MELE strongly suggests tail independence, since it hardly shows any increase away from the endpoints 0 and $\pi/2$. The conclusion is that extremely tall women are not extremely heavy and conversely. This nicely illustrates the difference between (in)dependence and tail (in)dependence. Similar results were obtained for the data for males.

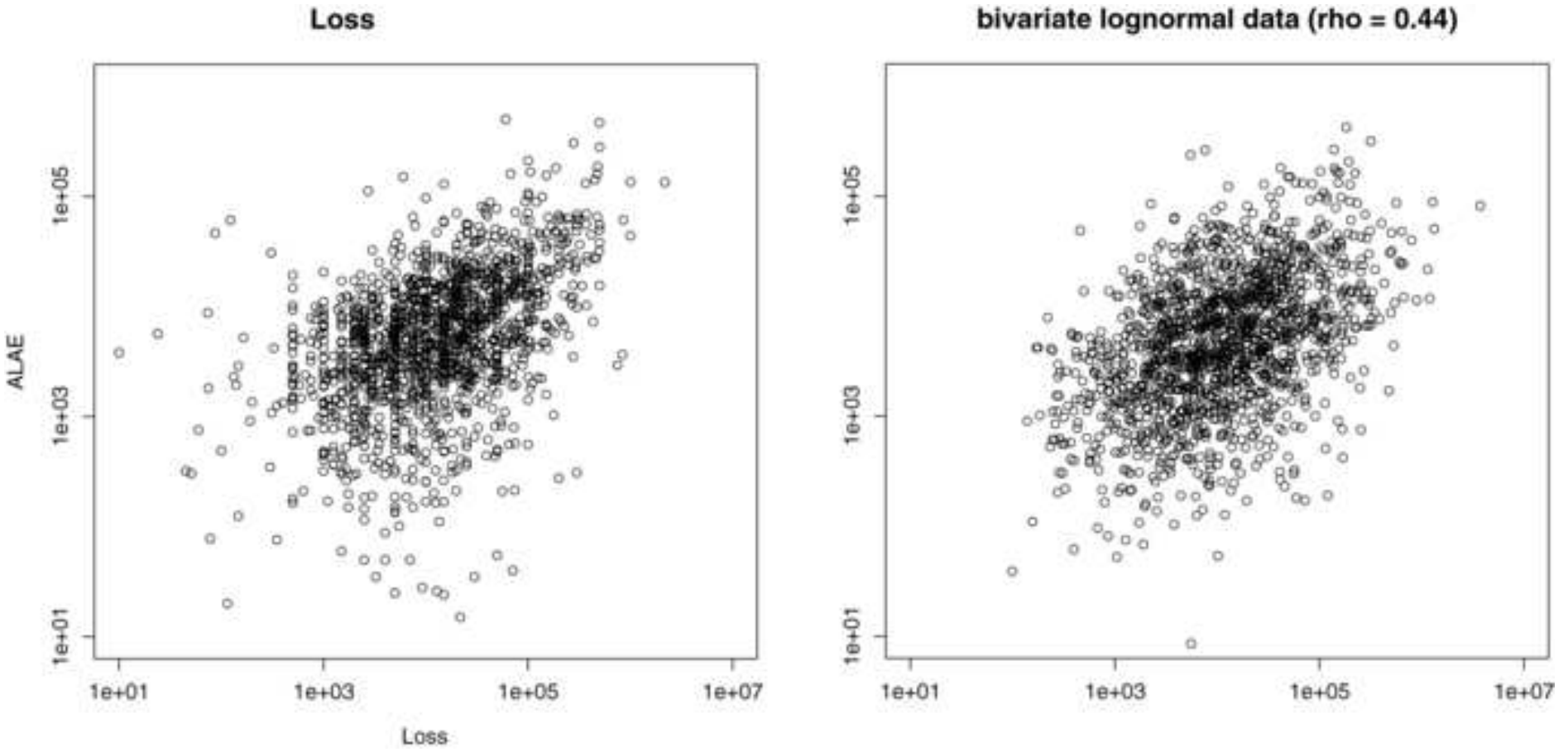

FIG. 6. *Left: Scatterplot of the Loss-ALAE data (axes on log-scale). Right: Scatterplot of a bivariate lognormal sample of the same size and which on a log-scale has the same mean vector and covariance matrix as the Loss-ALAE data.*

## 6. Proofs of Theorems 3.1–3.3.

PROOF OF THEOREM 3.1. A. We first prove weak convergence of the process $\sqrt{k}(\hat{\Phi}_p - \Phi_p)$ in $D[0, \pi/4]$. More precisely, with $\Delta \in \{1, \frac{1}{2}, \frac{1}{3}, \ldots\}$, we will show that for probabilistically equivalent versions of the processes

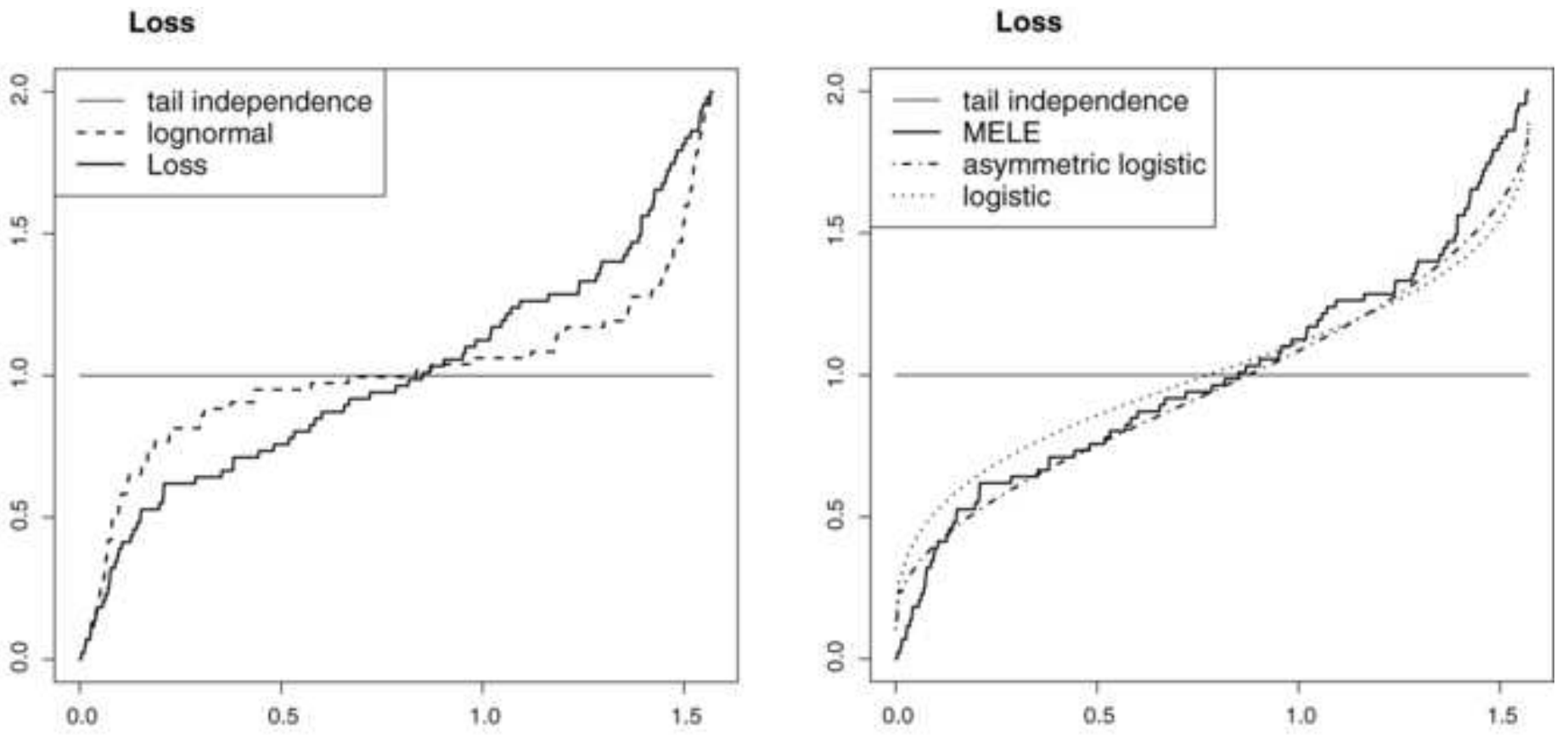

FIG. 7. *Left: MELE* ($p = 1$, $k = 40$) *of the spectral measure of the Loss-ALAE data and MELE of the bivariate lognormal data of Figure 6. Right: MELE of the spectral measure* ($p = 1$, $k = 40$) *of the Loss-ALAE data compared to the spectral measures of the fitted logistic and the fitted asymmetric logistic models.*

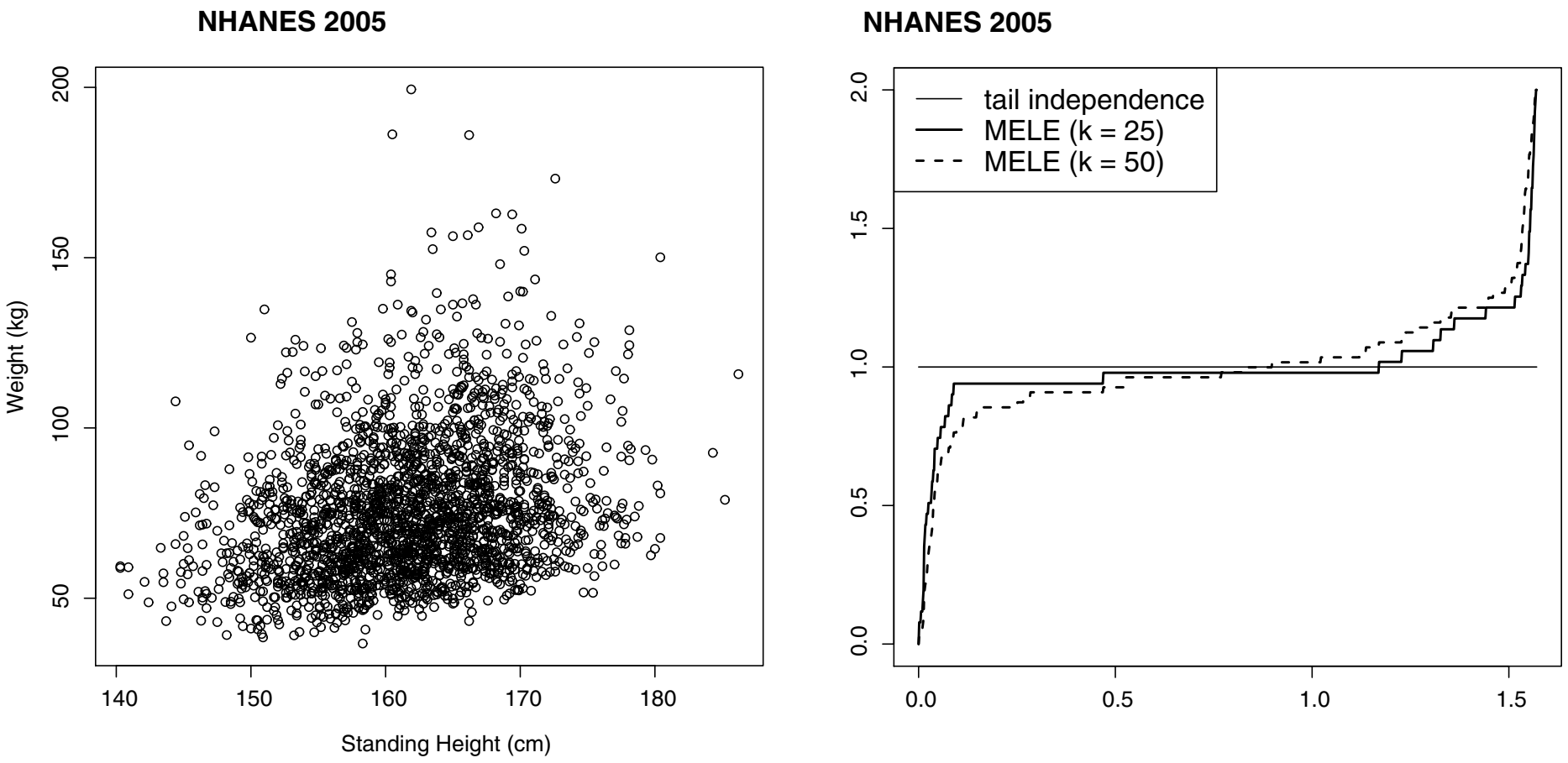

FIG. 8. *Scatterplot (left) and spectral measure estimates for $p=1$ (right) of the variables "standing height" and "weight" in the NHANES 2005-2006 body measurement data for females.*

involved and any $\varepsilon > 0$,

$$\lim_{\Delta\downarrow 0} \limsup_{n\to\infty} \Pr\Big\{ \sup_{\theta\in[0,\pi/4]} |\sqrt{k}\{\hat{\Phi}_p(\theta) - \Phi_p(\theta)\} - \{W_\Lambda(C_{p,\theta}) + Z_p(\theta)\}| \geq 3\varepsilon \Big\} = 0, \tag{6.1}$$

where $\hat{\Phi}_p = \hat{\Phi}_{p,\Delta}$ and $W_\Lambda = W_{\Lambda,\Delta}$. In part B below, we will prove weak convergence in $D[0,\pi/2]$.

Fix $\Delta \in \{1, \frac{1}{2}, \frac{1}{3}, \ldots\}$ and $M > 1$; later on, $M$ will be taken large. Let $\mathcal{A}' = \{A \cap A' : A, A' \in \mathcal{A}\}$, where $\mathcal{A} = \mathcal{A}(\Delta, M)$ is a Vapnik–Červonenkis (VC) class of sets defined as follows. For $m = 0, 1, 2, \ldots, \frac{1}{\Delta} - 1$ define

$$I_\Delta(m,\theta) = \begin{cases} [m\Delta x_p(\theta), (m+1)\Delta x_p(\theta)], & \text{if } \theta \in (0,\pi/4], \\ [0,\infty), & \text{if } \theta = 0 \text{ and } m = 0, \\ \varnothing, & \text{if } \theta = 0 \text{ and } m > 0; \end{cases}$$

$$J_\Delta(m) = [y_p(1 + (2^{1/p} - 1)(m+1)\Delta), y_p(1 + (2^{1/p} - 1)m\Delta)].$$

Set $\tilde{\mathcal{A}}$ to be the class containing all the following sets:

$$\bigcup_{m=0}^{1/\Delta - 1} \{(x,y) : x \in I_\Delta(m), 0 \leq y \leq x\tan\theta + B_m(x\tan\theta)^{1/16}\},$$
$$\text{for some } \theta \in [0,\pi/4] \text{ and } B_0, B_1, \ldots, B_{1/\Delta-1} \in [-1,1],$$

$$\bigcup_{m=0}^{1/\Delta - 1} \{(x,y) : x \in J_\Delta(m), x \geq x_p(\theta), 0 \leq y \leq y_p(x)(1 + K_m)\},$$

$$\text{for some } \theta \in [0, \pi/4] \text{ and } K_0, K_1, \dots, K_{1/\Delta - 1} \in \left[-\frac{1}{2}, \frac{1}{2}\right],$$

$$\{(x,y) : x \le a\}, \qquad \{(x,y) : y \le a\}$$

and

$$\{(x,y) : x \le a \text{ or } y \le a\}, \qquad \text{for some } a \in [0, M].$$

Next define $\tilde{\mathcal{A}}_s = \{A_s : A \in \tilde{\mathcal{A}}\}$, where, for $A \in \tilde{\mathcal{A}}$, $A_s = \{(x,y) : (y,x) \in A\}$. Finally define $\mathcal{A} = \tilde{\mathcal{A}} \cup \tilde{\mathcal{A}}_s$.

From $\lim_{t \downarrow 0} D_{1/t}(t) = 0$, $t^{-1} P([0,\infty] \times [0,t]) = \Lambda([0,\infty] \times [0,1]) = 1$ $(0 < t \le 1)$, and the homogeneity of $\lambda$ we obtain

$$\lim_{t \downarrow 0} \sup_{A \in \mathcal{A}'} |t^{-1} P(tA) - \Lambda(A)| = 0 \tag{6.2}$$

for all $\Delta \in \{1, \frac{1}{2}, \frac{1}{3}, \dots\}$ and $M > 1$. Theorem 3.1 in Einmahl (1997) now yields our basic convergence result: for a special construction (but keeping the same notation) we have

$$\sup_{A \in \mathcal{A}} \left| \sqrt{k} \left\{ \frac{n}{k} P_n\left(\frac{k}{n} A\right) - \frac{n}{k} P\left(\frac{k}{n} A\right) \right\} - W_\Lambda(A) \right| \xrightarrow{\text{a.s.}} 0, \qquad n \to \infty. \tag{6.3}$$

Throughout, we will work within this special construction.

In the sequel we can and will redefine $\hat{C}_{p,\theta}$, $\theta \in [0, \pi/4]$, by

$$\left\{ (x,y) : 0 \le x \le \infty, y \ge 0, y \le \frac{n}{k} Q_{2n}\left( (\tan\theta) \Gamma_{1n}\left(\frac{k}{n} x\right) \right), \right.$$
$$\left. y \le \frac{n}{k} Q_{2n}\left( \frac{k}{n} y_p\left( \frac{n}{k} \Gamma_{1n}\left(\frac{k}{n} x\right) \right) \right) \right\},$$

where $Q_{jn}$ is the quantile function corresponding to $\Gamma_{jn}$, $j = 1, 2$, with $Q_{jn}(y) := 0$ for $0 \le y \le (2n)^{-1}$ by convention. Define the marginal tail empirical processes by

$$w_{jn}(x) = \sqrt{k} \left\{ \frac{n}{k} \Gamma_{jn}\left(\frac{k}{n} x\right) - x \right\}, \qquad x \ge 0, \ j = 1, 2,$$

and the marginal tail quantile processes by

$$v_{jn}(x) = \sqrt{k} \left\{ \frac{n}{k} Q_{jn}\left(\frac{k}{n} x\right) - x \right\}, \qquad x \ge 0, \ j = 1, 2.$$

Note that $w_{jn}$ and $v_{jn}$ converge almost surely to $W_j$ and $-W_j$, respectively, for $j = 1, 2$, uniformly on $[0, M]$. Observe that for $x \ge 0$,

$$\frac{n}{k} Q_{2n}\left( (\tan\theta) \Gamma_{1n}\left(\frac{k}{n} x\right) \right) = x \tan\theta + \frac{z_{n,\theta}(x)}{\sqrt{k}}, \tag{6.4}$$

where

$$z_{n,\theta}(x) = w_{1n}(x)\tan\theta + v_{2n}\left(x\tan\theta + \frac{w_{1n}(x)}{\sqrt{k}}\tan\theta\right). \tag{6.5}$$

Also

$$\begin{aligned} &\frac{n}{k}Q_{2n}\left(\frac{k}{n}y_p\left(\frac{n}{k}\Gamma_{1n}\left(\frac{k}{n}x\right)\right)\right) \\ &\qquad = y_p\left(x + \frac{w_{1n}(x)}{\sqrt{k}}\right) + \frac{1}{\sqrt{k}}v_{2n}\left(y_p\left(x + \frac{w_{1n}(x)}{\sqrt{k}}\right)\right). \end{aligned} \tag{6.6}$$

We will treat the terms $V_{n,p}(\theta)$, $Y_{n,p}(\theta)$ and $r_{n,p}(\theta)$ from (3.4) in paragraphs A.1–A.3, respectively.

A.1. First we deal with $V_{n,p}(\theta)$ in (3.4). Set

$$\hat{C}_{p,\theta,1} = \{(x,y) \in \hat{C}_{p,\theta} : x < x_p(\theta)\} \quad \text{and} \quad \hat{C}_{p,\theta,2} = \hat{C}_{p,\theta} \setminus \hat{C}_{p,\theta,1}.$$

We focus on both sets separately when considering $V_{n,p}(\theta)$. For $p = \infty$, $\hat{C}_{p,\theta,1}$ has been dealt with in Einmahl, de Haan and Piterbarg (2001). We will omit the small modifications that are needed for general $p \in [1,\infty]$. However for $\hat{C}_{p,\theta,2}$, the case $p = \infty$ is trivial compared to $p \in [1,\infty)$. Therefore we will consider

$$V_{n,p,2}(\theta) := \sqrt{k}\left\{\frac{n}{k}P_n\left(\frac{k}{n}\hat{C}_{p,\theta,2}\right) - \frac{n}{k}P\left(\frac{k}{n}\hat{C}_{p,\theta,2}\right)\right\}$$

in detail now.

Recall $z_{n,\theta}(x)$ in (6.5) and define $s_n(x)$ through

$$y_p\left(x + \frac{w_{1n}(x)}{\sqrt{k}}\right) + \frac{1}{\sqrt{k}}v_{2n}\left(y_p\left(x + \frac{w_{1n}(x)}{\sqrt{k}}\right)\right) = y_p(x)\left(1 + \frac{s_n(x)}{\sqrt{k}}\right).$$

Further, put

$$W^+_{m,\Delta,\theta} = \sup_{\substack{x \in J_\Delta(m) \\ x \geq x_p(\theta)}}\left\{s_n(x) \wedge \left(\frac{z_{n,\theta}(x)}{y_p(x)} + \sqrt{k}\left(\frac{x\tan\theta}{y_p(x)} - 1\right)\right)\right\},$$

$$W^-_{m,\Delta,\theta} = \inf_{\substack{x \in J_\Delta(m) \\ x \geq x_p(\theta)}}\left\{s_n(x) \wedge \left(\frac{z_{n,\theta}(x)}{y_p(x)} + \sqrt{k}\left(\frac{x\tan\theta}{y_p(x)} - 1\right)\right)\right\}.$$

Set, for either choice of sign,

$$R^{\pm}_{m,\Delta,\theta} = \left\{(x,y) : x \in J_\Delta(m), x \geq x_p(\theta), 0 \leq y \leq y_p(x)\left(1 + \frac{W^{\pm}_{m,\Delta,\theta}}{\sqrt{k}}\right)\right\}$$

and

$$N^{\pm}_{\Delta,\theta} = \bigcup_{m=0}^{1/\Delta-1} R^{\pm}_{m,\Delta,\theta}.$$

We have

$$
\begin{aligned}
V_{n,p,2}(\theta) &\le \sqrt{k}\Big\{\frac{n}{k}P_n\Big(\frac{k}{n}N^+_{\Delta,\theta}\Big) - \frac{n}{k}P\Big(\frac{k}{n}N^+_{\Delta,\theta}\Big)\Big\} \\
&\quad + \sqrt{k}\frac{n}{k}P\Big(\frac{k}{n}(N^+_{\Delta,\theta} \setminus N^-_{\Delta,\theta})\Big) \qquad (6.7)\\
&=: V^+_{n,p,2}(\theta) + r_{n,p,2}(\theta);
\end{aligned}
$$

similarly

$$
\begin{aligned}
V_{n,p,2}(\theta) &\ge \sqrt{k}\Big\{\frac{n}{k}P_n\Big(\frac{k}{n}N^-_{\Delta,\theta}\Big) - \frac{n}{k}P\Big(\frac{k}{n}N^-_{\Delta,\theta}\Big)\Big\} \\
&\quad - \sqrt{k}\frac{n}{k}P\Big(\frac{k}{n}(N^+_{\Delta,\theta} \setminus N^-_{\Delta,\theta})\Big) \qquad (6.8)\\
&=: V^-_{n,p,2}(\theta) - r_{n,p,2}(\theta).
\end{aligned}
$$

We first deal with $r_{n,p,2}(\theta)$ and next with $V^{\pm}_{n,p,2}(\theta)$. Using (3.6) and well-known results on tail empirical and tail quantile processes [Csörgő and Horváth (1993) and Einmahl (1997)] we can show that, as $n \to \infty$,

$$
\sup_{\theta\in[0,\pi/4]} |r_{n,p,2}(\theta) - \sqrt{k}\Lambda(N^+_{\Delta,\theta} \setminus N^-_{\Delta,\theta})| \xrightarrow{p} 0. \qquad (6.9)
$$

Now consider

$$
\sup_{\theta\in[0,\pi/4]} \sqrt{k}\Lambda(N^+_{\Delta,\theta} \setminus N^-_{\Delta,\theta}).
$$

Set $c_m = y_p(1 + (2^{1/p} - 1)m\Delta) \vee x_p(\theta)$ and note that

$$
\sqrt{k}\Lambda(N^+_{\Delta,\theta} \setminus N^-_{\Delta,\theta}) \le \sqrt{k} \sum_{m=0}^{1/\Delta-1} \int_{c_{m+1}}^{c_m} \int_{y_p(x)(1+W^-_{m,\Delta,\theta}/\sqrt{k})}^{y_p(x)(1+W^+_{m,\Delta,\theta}/\sqrt{k})} \lambda(x,y)\,dy\,dx.
$$

Setting $y = y_p(x)(1 + z/\sqrt{k})$, we can rewrite the right-hand side of the previous display as

$$
\begin{aligned}
&\sum_{m=0}^{1/\Delta-1} \int_{c_{m+1}}^{c_m} y_p(x) \int_{W^-_{m,\Delta,\theta}}^{W^+_{m,\Delta,\theta}} \lambda\Big(x, y_p(x)\Big(1 + \frac{z}{\sqrt{k}}\Big)\Big)\,dz\,dx \\
&\quad = \sum_{m=0}^{1/\Delta-1} \int_{c_{m+1}}^{c_m} \int_{W^-_{m,\Delta,\theta}}^{W^+_{m,\Delta,\theta}} \frac{1}{1+z/\sqrt{k}} \lambda\Big(\frac{(x^p-1)^{1/p}}{1+z/\sqrt{k}}, 1\Big)\,dz\,dx \\
&\quad = \sum_{m=0}^{1/\Delta-1} \int_{W^-_{m,\Delta,\theta}}^{W^+_{m,\Delta,\theta}} \int_{c_{m+1}}^{c_m} \frac{1}{1+z/\sqrt{k}} \lambda\Big(\frac{(x^p-1)^{1/p}}{1+z/\sqrt{k}}, 1\Big)\,dx\,dz
\end{aligned}
$$

$$= \sum_{m=0}^{1/\Delta-1} \int_{W^-_{m,\Delta,\theta}}^{W^+_{m,\Delta,\theta}} \int_{(c^p_{m+1}-1)^{1/p}/(1+W^+_{m,\Delta,\theta}/\sqrt{k})}^{(c^p_m-1)^{1/p}/(1+W^-_{m,\Delta,\theta}/\sqrt{k})} \left( \frac{\{(1+z/\sqrt{k})v\}^p}{1+\{(1+z/\sqrt{k})v\}^p} \right)^{1-1/p} \times \lambda(v,1)\,dv\,dz.$$

The integrand is bounded by $\lambda(v,1)$, whence

$$\sqrt{k}\Lambda(N^+_{\Delta,\theta} \setminus N^-_{\Delta,\theta}) \le \max_{m\in\{0,1,\ldots,1/\Delta-1\}} (W^+_{m,\Delta,\theta} - W^-_{m,\Delta,\theta}) \times \sum_{m=0}^{1/\Delta-1} \int_{(c^p_{m+1}-1)^{1/p}/(1+W^+_{m,\Delta,\theta}/\sqrt{k})}^{(c^p_m-1)^{1/p}/(1+W^-_{m,\Delta,\theta}/\sqrt{k})} \lambda(v,1)\,dv. \tag{6.10}$$

We have

$$s_n(x) = \frac{1}{y_p(x)}\sqrt{k}\left\{ y_p\left(x + \frac{w_{1n}(x)}{\sqrt{k}}\right) - y_p(x) \right\} + \frac{1}{y_p(x)} v_{2n}\left( y_p\left( x + \frac{w_{1n}(x)}{\sqrt{k}} \right) \right). \tag{6.11}$$

Now from the behavior of tail empirical and tail quantile processes it readily follows that $\sup_{x\in[2^{1/p},\infty)} |s_n(x)| = O_p(1)$. Hence the right-hand side of (6.10) can be bounded, with probability tending to one, by

$$3 \max_{m\in\{0,1,\ldots,1/\Delta-1\}} (W^+_{m,\Delta,\theta} - W^-_{m,\Delta,\theta}) \int_0^\infty \lambda(v,1)\,dv.$$

As $\Lambda$ has uniform marginals, necessarily $\int_0^\infty \lambda(v,1)\,dv \le 1$. So in summary we have for fixed $\Delta$ and with probability tending to one,

$$\sup_{\theta\in[0,\pi/4]} \sqrt{k}\Lambda(N^+_{\Delta,\theta} \setminus N^-_{\Delta,\theta}) \le 3 \sup_{\theta\in[0,\pi/4]} \max_{m\in\{0,1,\ldots,1/\Delta-1\}} (W^+_{m,\Delta,\theta} - W^-_{m,\Delta,\theta}). \tag{6.12}$$

From the behavior of $s_n$ it follows that for any $\delta > 0$,

$$\lim_{\Delta\downarrow 0} \limsup_{n\to\infty} \Pr\left\{ \sup_{\theta\in[0,\pi/4]} \max_{m\in\{0,1,\ldots,1/\Delta-1\}} (W^+_{m,\Delta,\theta} - W^-_{m,\Delta,\theta}) \ge \delta \right\} = 0,$$

and hence, by (6.9),

$$\lim_{\Delta\downarrow 0} \limsup_{n\to\infty} \Pr\left\{ \sup_{\theta\in[0,\pi/4]} r_{n,p,2}(\theta) \ge \frac{\varepsilon}{2} \right\} = 0. \tag{6.13}$$

Next consider $V^\pm_{n,p,2}(\theta)$, for either choice of sign. Since

$$\lim_{n\to\infty} \Pr\{N^\pm_{\Delta,\theta} \in \tilde{\mathcal{A}}, \text{ for all } \theta \in [0,\pi/4]\} = 1,$$

we have, using (6.3),

$$\sup_{\theta\in[0,\pi/4]} |V^{\pm}_{n,p,2}(\theta) - W_\Lambda(N^{\pm}_{\Delta,\theta})| \xrightarrow{p} 0, \qquad n\to\infty. \tag{6.14}$$

But with similar calculations as for (6.12) we obtain

$$\Lambda(N^{\pm}_{\Delta,\theta} \triangle C_{p,\theta,2}) \le \frac{3}{\sqrt{k}} \max_{m\in\{0,1,\dots,1/\Delta-1\}} |W^{\pm}_{m,\Delta,\theta}|$$

with $C_{p,\theta,2} = \{(x,y)\in C_{p,\theta} : x \ge x_p(\theta)\}$. Since

$$\sup_{\theta\in[0,\pi/4]} \max_{m\in\{0,1,\dots,1/\Delta-1\}} |W^{\pm}_{m,\Delta,\theta}| = O_p(1), \qquad n\to\infty,$$

we have for any $\Delta\in\{1,\frac{1}{2},\frac{1}{3},\dots\}$,

$$\sup_{\theta\in[0,\pi/4]} \Lambda(N^{\pm}_{\Delta,\theta} \triangle C_{p,\theta,2}) \xrightarrow{p} 0, \qquad n\to\infty.$$

Hence, since $W_\Lambda$ is uniformly continuous on $\mathcal{A}$ with respect to the pseudo-metric $\Lambda(A\triangle A')$ for $A, A'\in\mathcal{A}$,

$$\sup_{\theta\in[0,\pi/4]} |W_\Lambda(N^{\pm}_{\Delta,\theta}) - W_\Lambda(C_{p,\theta,2})| \xrightarrow{p} 0, \qquad n\to\infty. \tag{6.15}$$

Combining (6.7), (6.8), (6.13), (6.14) and (6.15), we now have proven

$$\lim_{\Delta\downarrow 0}\limsup_{n\to\infty} \Pr\Big\{\sup_{\theta\in[0,\pi/4]} |V_{n,p,2}(\theta) - W_\Lambda(C_{p,\theta,2})| \ge \varepsilon\Big\} = 0.$$

This, in conjunction with the aforementioned result for $\hat{C}_{p,\theta,1}$, yields

$$\lim_{\Delta\downarrow 0}\limsup_{n\to\infty} \Pr\Big\{\sup_{\theta\in[0,\pi/4]} |V_{n,p}(\theta) - W_\Lambda(C_{p,\theta})| \ge 2\varepsilon\Big\} = 0. \tag{6.16}$$

A.2. Next we consider $Y_{n,p}(\theta) = \sqrt{k}\{\Lambda(\hat{C}_{p,\theta}) - \Lambda(C_{p,\theta})\}$. We will show that

$$\sup_{\theta\in[0,\pi/4]} |Y_{n,p}(\theta) - Z_p(\theta)| \xrightarrow{p} 0, \qquad n\to\infty. \tag{6.17}$$

Again, we will only consider $\hat{C}_{p,\theta,2}$. The other part, $\hat{C}_{p,\theta,1}$, can be handled as in Einmahl, de Haan and Piterbarg (2001); only minor modifications are needed.

So we will need to show that, as $n\to\infty$,

$$\begin{aligned}\sup_{\theta\in[0,\pi/4]} \bigg| &\sqrt{k}\{\Lambda(\hat{C}_{p,\theta,2}) - \Lambda(C_{p,\theta,2})\} \\ &- \int_{x_p(\theta)}^{\infty} \lambda(x, y_p(x))\{W_1(x)y_p'(x) - W_2(y_p(x))\}\,dx\bigg| \xrightarrow{p} 0.\end{aligned} \tag{6.18}$$

Observe, with $z_{n,\theta}$ and $s_n$ as in (6.5) and (6.11), respectively, that

$$\begin{aligned}(6.19)\qquad &\sqrt{k}\{\Lambda(\hat{C}_{p,\theta,2}) - \Lambda(C_{p,\theta,2})\}\\ &\quad = \sqrt{k}\int_{x_p(\theta)}^{\infty}\int_{y_p(x)}^{y_p(x)\{1+\check{s}_n(x)/\sqrt{k}\}} \lambda(x,y)\,dy\,dx,\end{aligned}$$

where

$$\check{s}_n(x) = s_n(x) \wedge \left\{\frac{z_{n,\theta}(x)}{y_p(x)} + \sqrt{k}\left(\frac{x\tan\theta}{y_p(x)} - 1\right)\right\}.$$

Since for fixed $x > x_p(\theta)$ the expression $\sqrt{k}(\frac{x\tan\theta}{y_p(x)} - 1)$ tends to infinity, it follows that we can (and will) replace $\check{s}_n(x)$ by $s_n(x)$ in the integral on the right-hand side of (6.19). Write

$$s(x) = \frac{1}{y_p(x)}\{W_1(x)y_p'(x) - W_2(y_p(x))\}.$$

Now

$$\begin{aligned}&\sup_{\theta\in[0,\pi/4]}\left|\sqrt{k}\int_{x_p(\theta)}^{\infty}\int_{y_p(x)}^{y_p(x)(1+s_n(x)/\sqrt{k})} \lambda(x,y)\,dy\,dx\right.\\ &\qquad\qquad \left. - \int_{x_p(\theta)}^{\infty} y_p(x)\lambda(x,y_p(x))s(x)\,dx\right|\\ &\quad \le \sup_{\theta\in[0,\pi/4]}\left|\sqrt{k}\int_{x_p(\theta)}^{\infty}\int_{y_p(x)(1+s(x)/\sqrt{k})}^{y_p(x)(1+s_n(x)/\sqrt{k})} \lambda(x,y)\,dy\,dx\right.\\ &\qquad\qquad + \sqrt{k}\int_{x_p(\theta)}^{\infty}\int_{y_p(x)}^{y_p(x)(1+s(x)/\sqrt{k})} \lambda(x,y)\,dy\,dx\\ &\qquad\qquad\qquad \left. - \int_{x_p(\theta)}^{\infty} y_p(x)\lambda(x,y_p(x))s(x)\,dx\right|\\ &\quad \le \sup_{\theta\in[0,\pi/4]}\left|\int_{x_p(\theta)}^{\infty}\int_{s(x)}^{s_n(x)} y_p(x)\lambda\left(x, y_p(x)\left(1+\frac{z}{\sqrt{k}}\right)\right)dz\,dx\right|\\ &\qquad + \sup_{\theta\in[0,\pi/4]}\left|\int_{x_p(\theta)}^{\infty}\int_{0}^{s(x)} y_p(x)\left[\lambda\left(x, y_p(x)\left(1+\frac{z}{\sqrt{k}}\right)\right)\right.\right.\\ &\qquad\qquad\qquad\qquad \left.\left. - \lambda(x,y_p(x))\right]dz\,dx\right|\\ &\quad =: T_1 + T_2.\end{aligned}$$

Since $\lambda(v,1) = v^{-1}\lambda(1,1/v)$ and by continuity of $\lambda$ on $[0,\infty)^2 \setminus \{(0,0)\}$, we have $\lim_{v\to\infty}\lambda(v,1) = 0$ and thus $\sup_{v\ge 0}\lambda(v,1) < \infty$. For some (large)

$M > 2$

$$T_1 \le \sup_{\theta\in[0,\pi/4]} \left| \int_{x_p(\theta)}^{M\vee x_p(\theta)} \int_{s(x)}^{s_n(x)} \frac{1}{1+z/\sqrt{k}} \lambda\left( \frac{(x^p-1)^{1/p}}{1+z/\sqrt{k}}, 1 \right) dz\,dx \right|$$
$$+ \left| \int_M^\infty \int_{s(x)}^{s_n(x)} \lambda\left( (x^p-1)^{1/p}, 1+\frac{z}{\sqrt{k}} \right) dz\,dx \right| =: T_{1,1} + T_{1,2}.$$

We first show

$$\sup_{2^{1/p}\le x\le M} |s_n(x) - s(x)| \xrightarrow{p} 0, \qquad n\to\infty. \tag{6.20}$$

Define

$$\tilde{s}_n(x) = \frac{y_p'(x)}{y_p(x)} w_{1n}(x) + \frac{1}{y_p(x)} v_{2n}\left( y_p\left( x + \frac{w_{1n}(x)}{\sqrt{k}} \right) \right).$$

Then it follows from the mean-value theorem and the almost sure convergence of $w_{1n}$ to $W_1$, uniformly on $[0, M]$, that

$$\sup_{2^{1/p}\le x\le M} |s_n(x) - \tilde{s}_n(x)| \xrightarrow{p} 0, \qquad n\to\infty.$$

It also follows easily that

$$\sup_{2^{1/p}\le x\le M} |\tilde{s}_n(x) - s(x)| \xrightarrow{p} 0, \qquad n\to\infty,$$

whence (6.20). We have with probability tending to one,

$$T_{1,1} \le 2M \sup_{v\ge 0} \lambda(v,1) \sup_{2^{1/p}\le x\le M} |s_n(x) - s(x)|,$$

which, because of (6.20), tends to 0 in probability (for any $M > 2$). Let $\kappa > 0$ and set $\delta = \sqrt{\kappa}/2$. Using again (6.20) and the behavior of $W_1$ near infinity, we see that for large enough $M$ and with probability tending to one,

$$\begin{aligned}
T_{1,2} &\le \int_M^\infty \int_{-W_2(1)-\delta}^{-W_2(1)+\delta} \lambda\left( (x^p-1)^{1/p}, 1+\frac{z}{\sqrt{k}} \right) dz\,dx \\
&= \int_{-W_2(1)-\delta}^{-W_2(1)+\delta} \int_M^\infty \frac{1}{1+z/\sqrt{k}} \lambda\left( \frac{(x^p-1)^{1/p}}{1+z/\sqrt{k}}, 1 \right) dx\,dz \\
&\le \int_{-W_2(1)-\delta}^{-W_2(1)+\delta} \int_{1/2(M^p-1)^{1/p}}^\infty \left( \frac{\{(1+z/\sqrt{k})v\}^p}{1+\{(1+z/\sqrt{k})v\}^p} \right)^{1-1/p} \lambda(v,1)\,dv\,dz \\
&\le \int_{-W_2(1)-\delta}^{-W_2(1)+\delta} \delta\,dz = 2\delta^2 = \kappa/2,
\end{aligned}$$

whence

$$\lim_{M\to\infty} \limsup_{n\to\infty} \Pr\{T_1 \ge \kappa\} = 0. \tag{6.21}$$

Now consider $T_2$. Write $\|s\| = \sup_{2^{1/p} \le x \le \infty} |s(x)|$ and

$$D_n = \sup_{\substack{x \ge 2^{1/p} \\ -\|s\| \le z \le \|s\|}} \left| \lambda\left( (x^p - 1)^{1/p}, 1 + \frac{z}{\sqrt{k}} \right) - \lambda((x^p - 1)^{1/p}, 1) \right|.$$

For $M > 2^{1/p}$,

$$T_2 = \sup_{\theta \in [0,\pi/4]} \left| \int_{x_p(\theta)}^{\infty} \int_0^{s(x)} \left\{ \lambda\left( (x^p - 1)^{1/p}, 1 + \frac{z}{\sqrt{k}} \right) - \lambda((x^p - 1)^{1/p}, 1) \right\} dz\, dx \right|$$

$$\le \int_{2^{1/p}}^{M} \int_{-\|s\|}^{\|s\|} D_n \, dz\, dx$$

$$+ \int_M^{\infty} \int_{-\|s\|}^{\|s\|} \left\{ \lambda\left( (x^p - 1)^{1/p}, 1 + \frac{z}{\sqrt{k}} \right) + \lambda((x^p - 1)^{1/p}, 1) \right\} dz\, dx.$$

Clearly, as $n \to \infty$,

$$\int_{2^{1/p}}^{M} \int_{-\|s\|}^{\|s\|} D_n \, dz\, dx \le 2M\|s\| D_n \xrightarrow{p} 0$$

and also, with probability tending to one,

$$\int_M^{\infty} \int_{-\|s\|}^{\|s\|} \left\{ \lambda\left( (x^p - 1)^{1/p}, 1 + \frac{z}{\sqrt{k}} \right) + \lambda((x^p - 1)^{1/p}, 1) \right\} dz\, dx$$

$$= \int_M^{\infty} \int_{-\|s\|}^{\|s\|} \frac{1}{1 + z/\sqrt{k}} \lambda\left( \frac{(x^p - 1)^{1/p}}{1 + z/\sqrt{k}}, 1 \right) dz\, dx$$

$$+ \int_M^{\infty} \int_{-\|s\|}^{\|s\|} \lambda((x^p - 1)^{1/p}, 1)\, dz\, dx$$

$$\le \int_{-\|s\|}^{\|s\|} \int_{1/2(M^p - 1)^{1/p}}^{\infty} \lambda(v, 1)\, dv\, dz + \int_{-\|s\|}^{\|s\|} \int_{(M^p - 1)^{1/p}}^{\infty} \lambda(u, 1)\, du\, dz$$

$$\le 4\|s\| \int_{1/2(M^p - 1)^{1/p}}^{\infty} \lambda(v, 1)\, dv.$$

As a result,

$$\lim_{M \to \infty} \limsup_{n \to \infty} \Pr\{T_2 \ge \kappa\} = 0. \tag{6.22}$$

Combining (6.21) and (6.22) yields (6.18), which, in conjunction with the aforementioned result for $\Lambda(\hat{C}_{p,\theta,1})$, yields (6.17).

A.3. We now consider $r_{n,p}(\theta)$ in (3.4). From (6.4), (6.6), (3.6) and the behavior of tail empirical and tail quantile processes, it follows that

$$\sup_{\theta \in [0,\pi/4]} |r_{n,p}(\theta)| \xrightarrow{p} 0 \qquad \text{as } n \to \infty. \tag{6.23}$$

Combining (6.16), (6.17) and (6.23) yields (6.1). So actually we proved the theorem for $\theta \in [0, \pi/4]$.

B. Observe that, using a symmetry argument, it rather easily follows from (6.1) with $\theta = \pi/4$ that

$$\lim_{\Delta\downarrow 0} \limsup_{n\to\infty} \Pr[|\sqrt{k}\{\hat{\Phi}_p(\pi/2) - \Phi_p(\pi/2)\} - \{W_\Lambda(C_{p,\pi/2}) + Z_p(\pi/2)\}| \geq 6\varepsilon] = 0.$$

Observe in particular that the first term of $Z_{c,p}(\pi/4)$ cancels out with the similar term coming from the mirror image (with respect to the line $y = x$) of $C_{p,\pi/4}$. By a similar symmetry argument, observing that for $\theta \in (\pi/4, \pi/2)$ (the closure of) $C_{p,\pi/2} \setminus C_{p,\theta}$ is the mirror image of $C_{p,\pi/2-\theta}$, it follows that

$$\lim_{\Delta\downarrow 0} \limsup_{n\to\infty} \Pr\Big[\sup_{\theta\in(\pi/4,\pi/2]} |\sqrt{k}\{\hat{\Phi}_p(\theta) - \Phi_p(\theta)\} - \{W_\Lambda(C_{p,\theta}) + Z_p(\theta)\}| \geq 9\varepsilon\Big] = 0. \tag{6.24}$$

Combining (6.1) and (6.24) completes the proof. □

PROOF OF THEOREM 3.2. The proof of this theorem follows in the same way as that of Theorem 3.1; only small adaptations are needed, including the obvious adaptation of the VC class $\mathcal{A}$. The main difference between both results is the weaker condition (3.8) which allows $\Lambda$ to put mass on $\{\infty\} \times [0, \infty)$ or $[0, \infty) \times \{\infty\}$; on the other hand $\theta$ is bounded away from 0 and $\pi/2$ in the present result. In the limit process, the term $W_\Lambda(C_{p,\theta})$ stays the same as in Theorem 3.1 but with weaker conditions on $\Lambda$; the term $Z_p(\theta) = Z_{c,p}(\theta) + Z_d(\theta)$ may now be different from that in Theorem 3.1, since there $Z_d = 0$, which might not be the case here. Therefore, we confine ourselves to explaining how condition (3.8) is set to use and to the adaptation of that part of the proof that deals with $Z_d$.

Condition (3.8) implies that for some sequence $T_n$

$$\sqrt{k} D_{T_n}(k/n) + \sqrt{k}/T_n + 1/T_n^{1/2} \to 0, \qquad n \to \infty. \tag{6.25}$$

We focus on the bias term $\sup_{\theta\in[\eta,\pi/4]} |r_{n,p}(\theta)|$; see (3.4). For $\theta \in [\eta, \pi/4]$, write $\hat{C}_{p,\theta} = C_1 \cup C_2 \cup C_3$, where

$$C_1 = \hat{C}_{p,\theta} \cap ([0, T_n] \times [0, \infty]),$$

$$C_2 = \hat{C}_{p,\theta} \cap \left([T_n, \infty] \times \left[0, \frac{n}{k} Q_{2n}\left(\frac{k}{n}\right)\right]\right),$$

$$C_3 = \hat{C}_{p,\theta} \cap \left([T_n, \infty] \times \left[\frac{n}{k} Q_{2n}\left(\frac{k}{n}\right), \infty\right)\right).$$

By the triangle inequality the bias term can be split up into three terms, based on $C_1$, $C_2$ and $C_3$, respectively. The first one of these terms converges to zero in probability, because the first term in (6.25) tends to 0. Using

$$\frac{n}{k}P\left([0,\infty]\times\left[0,Q_{2n}\left(\frac{k}{n}\right)\right]\right)=\frac{n}{k}Q_{2n}\left(\frac{k}{n}\right)=\Lambda\left([0,\infty]\times\left[0,\frac{n}{k}Q_{2n}\left(\frac{k}{n}\right)\right]\right),$$

the second one can be handled similarly. For the third term we replace the difference in the definition of $r_{n,p}(\theta)$ by a sum and deal with both terms obtained from this sum separately. Using the behavior of tail empirical and tail quantile processes we obtain the convergence of both these terms from the convergence to 0 of the second and third term in (6.25).

Recall that $Z_d(\theta)=-\Phi_p(\{0\})W_2(1)$. We have to show the following analogue of (6.18):

$$\sup_{\theta\in[\eta,\pi/4]}\biggl|\sqrt{k}\{\Lambda(\hat{C}_{p,\theta,2})-\Lambda(C_{p,\theta,2})\} - \int_{x_p(\theta)}^{\infty}\lambda(x,y_p(x))\{W_1(x)y_p'(x)-W_2(y_p(x))\}\,dx + \Phi_p(\{0\})W_2(1)\biggr| \xrightarrow{p} 0. \tag{6.26}$$

In view of the proof of (6.18), the proof of (6.26) is complete if we show that, as $n\to\infty$,

$$\sup_{\theta\in[\eta,\pi/4]}|\Phi_p(\{0\})v_{2n}(1)+\Phi_p(\{0\})W_2(1)|=\Phi_p(\{0\})|v_{2n}(1)+W_2(1)|\xrightarrow{p}0.$$

But this immediately follows from (6.3). □

Proof of Theorem 3.3. The proof of Theorem 3.3 goes along the same lines of those of Theorems 3.1–3.2. Observe that we only have to consider the process $\sqrt{k}(\hat{\Phi}_p-\Phi_p)$ on $[\eta_n,\pi/2-\eta_n]$ and at $\pi/2$, since on $[0,\eta_n)$ and $(\pi/2-\eta_n,\pi/2)$ the process is constant and the limit process is continuous on $[0,\pi/2)$. Then we are in a similar situation as in Theorem 3.2, but now the interval under consideration depends on $n$ and converges to $(0,\pi/2)$.

The essential difference lies in the VC class $\mathcal{A}$. If we would adapt the VC class in the proof of Theorem 3.1 in the obvious way, that is, restrict $\theta$ to $[\eta_n,\pi/2-\eta_n]$, the VC class would depend on $n$ and hence Theorem 3.1 in Einmahl (1997) would not be applicable. We will, however, consider the VC class that is obtained from $\mathcal{A}$ of our Theorem 3.1 by omitting $\theta=0$. Of course, (6.2) does not necessarily hold for this new class, but it can be shown to hold when we replace $\frac{n}{k}P(\frac{k}{n}\cdot)$ by $\widetilde{P}_{(n)}$, the measure that is obtained from $\frac{n}{k}P(\frac{k}{n}\cdot)$ by projecting the probability mass of

$$\frac{k}{n}([T_n, n/k] \times ([0, 1-k^{-1/4}] \cup [1+k^{-1/4}, 3])) \tag{6.27}$$

on the axis $\{\infty\} \times [0,\infty)$, and by projecting the probability mass of

$$\frac{k}{n}(([0, 1-k^{-1/4}] \cup [1+k^{-1/4}, 3]) \times [T_n, n/k]) \tag{6.28}$$

on the axis $[0,\infty) \times \{\infty\}$; here $T_n \geq 2/\eta_n$ is a sequence of $T$s for which (3.11) holds. The points $\frac{n}{k}(U_{i1}, U_{i2})$, $i = 1, \ldots, n$, in the region (6.27) or (6.28) are projected on $\{\infty\} \times [0,\infty)$ or $[0,\infty) \times \{\infty\}$ in a similar way, that is, are replaced by $(\infty, \frac{n}{k}U_{i2})$ or $(\frac{n}{k}U_{i1}, \infty)$, respectively. It is easily seen that, with probability tending to one, this projection does not change the processes involved in the result. □

## 7. Proofs of Theorems 4.1–4.2.

Proof of Theorem 4.1. Equation (4.9) is an immediate consequence of (4.8), so we focus on (4.8).

Similarly it is immediate from Theorem 3.1 that, in $D[0, \pi/2]$ and as $n \to \infty$,

$$\sqrt{k}(\hat{Q}_p - Q_p) \xrightarrow{d} \frac{\Phi_p(\pi/2)\alpha_p - \alpha_p(\pi/2)\Phi_p}{\Phi_p^2(\pi/2)} = \beta_p. \tag{7.1}$$

Recall the definition of $f$ in (4.2) and observe that $\sup_{0\leq\theta\leq\pi/2} |f(\theta)| = 1$. Put $A_{in} = f(\Theta_{in})$ for $i \in I_n$. By (7.1) and since $Q_p(\{\pi/4\}) < 1$, necessarily $\Pr[\exists i \in I_n : A_{in} \neq 0] \to 1$ as $n \to \infty$.

Consider the random function

$$\Psi_n(\mu) = \frac{1}{N_n} \sum_{i\in I_n} \frac{A_{in}}{1+\mu A_{in}}, \qquad -1 < \mu < 1.$$

The derivative of $\Psi_n$ is

$$\Psi_n'(\mu) = -\frac{1}{N_n} \sum_{i\in N_n} \frac{A_{in}^2}{(1+\mu A_{in})^2}.$$

Hence, on the event $\{\exists i \in I_n : A_{in} \neq 0\}$, the function $\Psi_n$ is strictly decreasing and there can be at most one $\tilde{\mu}_n \in (-1, 1)$ with $\Psi_n(\tilde{\mu}_n) = 0$.

If $g : [0, \pi/2] \to \mathbb{R}$ is absolutely continuous with Radon–Nikodym derivative $g'$, then by Fubini's theorem,

$$\frac{1}{N_n} \sum_{i\in I_n} g(\Theta_{in}) = \int_{[0,\pi/2]} g(\theta)\hat{Q}_p(d\theta)$$

$$= g(\pi/2) - \int_{[0,\pi/2]} \int_\theta^{\pi/2} g'(\vartheta)\, d\vartheta \hat{Q}_p(d\theta)$$

$$= g(\pi/2) - \int_0^{\pi/2} \hat{Q}_p(\vartheta) g'(\vartheta)\, d\vartheta.$$

Since similarly $\int g\, dQ_p = g(\pi/2) - \int_0^{\pi/2} Q_p(\vartheta) g'(\vartheta)\, d\vartheta$, by (7.1),

$$\begin{aligned} &\sqrt{k}\biggl(\frac{1}{N_n}\sum_{i\in I_n} g(\Theta_{in}) - \int_{[0,\pi/2]} g\, dQ_p\biggr) \\ &\quad = -\int_0^{\pi/2} \sqrt{k}\{\hat{Q}_p(\theta) - Q_p(\theta)\} g'(\theta)\, d\theta \qquad (7.2) \\ &\quad \xrightarrow{d} -\int_0^1 \beta_p(\theta) g'(\theta)\, d\theta, \qquad n\to\infty. \end{aligned}$$

Here we used the fact that the linear functional sending $x \in D[0,\pi/2]$ to $\int_0^{\pi/2} x(\theta) g'(\theta)\, d\theta$ is bounded.

Since $1/(1+x) = 1 - x + x^2/(1+x)$ for $x \neq -1$, we have

$$\begin{aligned} \Psi_n(\mu) &= \frac{1}{N_n}\sum_{i\in I_n} A_{in}\biggl(1 - \mu A_{in} + \frac{\mu^2 A_{in}^2}{1+\mu A_{in}}\biggr) \\ &= \frac{1}{N_n}\sum_{i\in I_n} A_{in} - \mu\frac{1}{N_n}\sum_{i\in I_n} A_{in}^2 + \mu^2\frac{1}{N_n}\sum_{i\in I_n}\frac{A_{in}^3}{1+\mu A_{in}}. \end{aligned}$$

Define

$$\bar{\mu}_n = \frac{1}{N_n}\sum_{i\in I_n} A_{in}\Big/\frac{1}{N_n}\sum_{i\in I_n} A_{in}^2 = \int f\, d\hat{Q}_p\Big/\int f^2\, d\hat{Q}_p.$$

Since $\int f\, dQ_p = 0$ and $\int f^2\, dQ_p > 0$, by (7.2), $\bar{\mu}_n = O_p(k^{-1/2})$ as $n\to\infty$. We have

$$\Psi_n(0) = \frac{1}{N_n}\sum_{i\in I_n} A_{in}$$

as well as

$$\begin{aligned} \Psi_n(2\bar{\mu}_n) &= -\frac{1}{N_n}\sum_{i\in I_n} A_{in} + 4\bar{\mu}_n^2\frac{1}{N_n}\sum_{i\in I_n}\frac{A_{in}^3}{1+2\bar{\mu}_n A_{in}} \\ &= -\frac{1}{N_n}\sum_{i\in I_n} A_{in}\cdot\biggl(1 - 4\bar{\mu}_n\frac{\sum_{i\in I_n} A_{in}^3/(1+2\bar{\mu}_n A_{in})}{\sum_{i\in I_n} A_{in}^2}\biggr). \end{aligned}$$

Because $\bar{\mu}_n = o_p(1)$, $|A_{in}| \le 1$ and $N_n^{-1}\sum_{i\in I_n} A_{in}^2 \xrightarrow{p} \int f^2\, dQ_p > 0$, we obtain

$$\lim_{n\to\infty} \Pr[|2\bar{\mu}_n| < 1, \Psi_n(0)\Psi_n(2\bar{\mu}_n) \le 0] = 1.$$

Since moreover, with probability tending to one, $\Psi_n$ is continuous and decreasing,

$$\lim_{n\to\infty} \Pr[\text{there exists a unique } \tilde{\mu}_n \in (-1,1) \text{ such that } \Psi_n(\tilde{\mu}_n) = 0] = 1.$$

Also, $\Pr(|\tilde{\mu}_n| \le 2|\bar{\mu}_n|) \to 1$ and thus $\tilde{\mu}_n = O_p(k^{-1/2})$ as $n \to \infty$. We have

$$0 = \Psi_n(\tilde{\mu}_n) = \frac{1}{N_n}\sum_{i\in I_n} A_{in} - \tilde{\mu}_n \frac{1}{N_n}\sum_{i\in I_n} A_{in}^2 + \tilde{\mu}_n^2 \frac{1}{N_n}\sum_{i\in I_n} \frac{A_{in}^3}{1+\tilde{\mu}_n A_{in}},$$

whence

$$\tilde{\mu}_n = \bar{\mu}_n + \tilde{\mu}_n^2 \frac{\sum_{i\in I_n} A_{in}^3/(1+\tilde{\mu}_n A_{in})}{\sum_{i\in I_n} A_{in}^2} = \bar{\mu}_n + O_p(k^{-1}), \qquad n \to \infty.$$

Define

$$\check{\mu}_n := \frac{1}{N_n}\sum_{i\in I_n} A_{in} \Big/ \int_{[0,\pi/2]} f^2(\theta) Q(d\theta).$$

Since $N_n^{-1}\sum_i A_{in} = O_p(k^{-1/2})$ and $N_n^{-1}\sum_i A_{in}^2 = \int f^2\, dQ_p + O_p(k^{-1/2})$, we have $\bar{\mu}_n = \check{\mu}_n + O_p(k^{-1})$ and thus

$$\tilde{\mu}_n = \check{\mu}_n + O_p(k^{-1}), \qquad n \to \infty.$$

Put

$$\check{p}_{in} := \frac{1}{N_n}(1 - \check{\mu}_n A_{in}), \qquad i \in I_n;$$

$$\check{Q}_p(\theta) := \sum_{i\in I_n} \check{p}_{in} \mathbf{1}(\Theta_{in} \le \theta), \qquad \theta \in [0,\pi/2].$$

Then

$$\begin{aligned}\tilde{Q}_p(\theta) - \check{Q}_p(\theta) &= \frac{1}{N_n}\sum_{i\in I_n}\left(\frac{1}{1+\tilde{\mu}_n A_{in}} - (1-\check{\mu}_n A_{in})\right)\mathbf{1}(\Theta_{in} \le \theta) \\ &= \frac{1}{N_n}\sum_{i\in I_n} \frac{(\check{\mu}_n - \tilde{\mu}_n)A_{in} + \tilde{\mu}_n\check{\mu}_n A_{in}^2}{1+\tilde{\mu}_n A_{in}} \mathbf{1}(\Theta_{in} \le \theta).\end{aligned}$$

Since both $\check{\mu}_n - \tilde{\mu}_n$ and $\tilde{\mu}_n\check{\mu}_n$ are $O_p(k^{-1})$,

$$\sup_{\theta\in[0,\pi/2]} |\tilde{Q}_p(\theta) - \check{Q}_p(\theta)| = O_p(k^{-1}), \qquad n \to \infty. \tag{7.3}$$

Therefore, as $n \to \infty$ and uniformly in $\theta \in [0,\pi/2]$,

$$\begin{aligned}&\sqrt{k}\{\tilde{Q}_p(\theta) - Q_p(\theta)\} \\ &\quad = \sqrt{k}\{\check{Q}_p(\theta) - Q_p(\theta)\} + O_p(k^{-1/2}) \\ &\quad = \sqrt{k}\{\hat{Q}_p(\theta) - Q_p(\theta)\} - \sqrt{k}\check{\mu}_n \hat{I}(\theta) + O_p(k^{-1/2}),\end{aligned} \tag{7.4}$$

where

$$\hat{I}(\theta) = \frac{1}{N_n} \sum_{i \in I_n} A_{in} \mathbf{1}(\Theta_{in} \le \theta), \qquad \theta \in [0, \pi/2].$$

The function $f$ is absolutely continuous; denote its Radon–Nikodym derivative by $f'$. By Fubini's theorem, for $\theta \in [0, \pi/2]$,

$$\hat{I}(\theta) = \int_{[0,\theta]} f \, d\hat{Q}_p = f(\theta)\hat{Q}_p(\theta) - \int_0^\theta \hat{Q}_p(\vartheta) f'(\vartheta) \, d\vartheta,$$

$$I(\theta) = \int_{[0,\theta]} f \, dQ_p = f(\theta)Q_p(\theta) - \int_0^\theta Q_p(\vartheta) f'(\vartheta) \, d\vartheta.$$

As a result,

$$\sup_{\theta \in [0,\pi/2]} |\hat{I}(\theta) - I(\theta)| = O_p(k^{-1/2}), \qquad n \to \infty. \tag{7.5}$$

Moreover, by (7.2) with $f = g$,

$$\sqrt{k}\check{\mu}_n = -\frac{1}{\int f^2 \, dQ_p} \int_0^{\pi/2} \sqrt{k}\{\hat{Q}_p(\theta) - Q_p(\theta)\} f'(\theta) \, d\theta. \tag{7.6}$$

Write $\beta_{n,p} = \sqrt{k}(\hat{Q}_p - Q_p)$. Combine (7.4), (7.5) and (7.6) to see that

$$\sqrt{k}(\tilde{Q}_p - Q_p) = \beta_{n,p} + \frac{\int \beta_{n,p} \, df}{\int f^2 \, dQ_p} I + O_p(k^{-1/2}) \tag{7.7}$$

as $n \to \infty$. Since the linear operator

$$D[0, \pi/2] \to D[0, \pi/2] : x \mapsto x + \frac{\int x \, df}{\int f^2 \, dQ_p} I$$

is bounded, (7.1) and (7.7) imply (4.8). $\square$

Proof of Theorem 4.2. It is immediate from Theorem 3.3 that, in $D[0, \pi/2]$ and as $n \to \infty$,

$$\sqrt{k}(\hat{Q}_p - Q_p) \circ \tau_n \xrightarrow{d} \beta_p.$$

Now the proof of Theorem 4.1 applies here as well, except for one change: we have to check that (7.2) still holds. But this follows from the fact that

$$\left| \int_0^{\pi/2} \sqrt{k}(\hat{Q}_p - Q_p)(\theta) g'(\theta) \, d\theta - \int_0^{\pi/2} \sqrt{k}(\hat{Q}_p - Q_p)(\tau_n(\theta)) g'(\theta) \, d\theta \right|$$

is bounded by

$$2\sqrt{k}\eta_n \sup_{\theta \in [0,\eta_n] \cup [\pi/2 - \eta_n, \pi/2]} |g'(\theta)|,$$

which by assumption tends to zero as $n \to \infty$ provided that $g'$ is bounded in the neighborhood of 0 and $\pi/2$. This is the case for $g = f$ and $g = f^2$, the only functions to which (7.2) is to be applied. $\square$

**Acknowledgment.** We are grateful to an editor and a referee for several thoughtful comments.

## REFERENCES

Abdous, B. and Ghoudi, K. (2005). Nonparametric estimators of multivariate extreme-dependence functions. *J. Nonparametr. Statist.* **17** 915–935. MR2192166

Beirlant, J., Goegebeur, Y., Segers, J. and Teugels, J. (2004). *Statistics of Extremes: Theory and Applications*. Wiley, Chichester. MR2108013

Capéraà, P. and Fougères, A.-L. (2000). Estimation of a bivariate extreme-value distribution. *Extremes* **3** 311–329. MR1870461

Coles, S. (2001). *An Introduction to Statistical Modelling of Extreme Values*. Springer, New York. MR1932132

Coles, S. and Tawn, J. (1991). Modelling extreme-multivariate events. *J. Roy. Statist. Soc. Ser. B* **53** 377–392. MR1108334

Csörgő, M. and Horváth, L. (1993). *Weighted Approximations in Probability and Statistics*. Wiley, Chichester. MR1215046

Drees, H. and Huang, X. (1998). Best attainable rates of convergence for estimates of the stable tail dependence function. *J. Multivariate Anal.* **64** 25–47. MR1619974

Einmahl, J. (1997). Poisson and Gaussian approximation of weighted local empirical processes. *Stochastic Process. Appl.* **70** 31–58. MR1472958

Einmahl, J., de Haan, L. and Li, D. (2006). Weighted approximations of tail copula processes with application to testing the bivariate extreme-value condition. *Ann. Statist.* **34** 1987–2014. MR2283724

Einmahl, J., de Haan, L. and Piterbarg, V. (2001). Nonparametric estimation of the spectral measure of an extreme-value distribution. *Ann. Statist.* **29** 1401–1423. MR1873336

Einmahl, J., de Haan, L. and Sinha, A. (1997). Estimating the spectral measure of an extreme-value distribution. *Stochastic Process. Appl.* **70** 143–171. MR1475660

Einmahl, J., Krajina, A. and Segers, J. (2008). A method of moments estimator for tail dependence. *Bernoulli* **14** 1003–1026.

Frees, E. W. and Valdez, E. A. (1998). Understanding relationships using copulas. *N. Am. Actuar. J.* **2** 1–25. MR1988432

Genest, C., Ghoudi, K. and Rivest, L.-P. (1998). Discussion of Frees and Valdez (1998). *N. Am. Actuar. J.* **2** 143–149. MR2011244

Gumbel, E. J. (1960). Bivariate exponential distributions. *J. Amer. Statist. Assoc.* **55** 698–707. MR0116403

de Haan, L. and Ferreira, A. (2006). *Extreme Value Theory: An Introduction*. Springer, New York. MR2234156

de Haan, L., Neves, C. and Peng, L. (2008). Parametric tail copula estimation and model testing. *J. Multivariate Anal.* **99** 1260–1275. MR2419346

de Haan, L. and Resnick, S. (1977). Limit theory for multidimensional sample extremes. *Z. Wahrsch. Verw. Gebiete* **40** 317–337. MR0478290

Huang, X. (1992). *Statistics of Bivariate Extremes*. Ph.D. thesis, Erasmus Univ. Rotterdam.

Joe, H., Smith, R. L. and Weissman, I. (1992). Bivariate threshold methods for extremes. *J. Roy. Statist. Soc. Ser. B* **54** 171–183. MR1157718

Ledford, A. W. and Tawn, J. A. (1996). Statistics for near independence in multivariate extreme values. *Biometrika* **83** 169–187. MR1399163

National Center for Health Statistics (2007). *National Health and Nutrition Examination Survey 2005–2006*. Available at http://www.cdc.gov/nchs/about/major/nhanes/nhanes2005-2006/exam05_06.htm.

Owen, A. (2001). *Empirical Likelihood*. Chapman & Hall/CRC, Boca Raton, FL.

Pickands, J. (1981). Multivariate extreme-value distributions. In *Proceedings of the 43rd session of the International Statistical Institute, Vol. 2 (Buenos Aires, 1981)* **49** 859–878, 894–902. With a discussion. MR0820979

Tawn, J. A. (1988). Bivariate extreme-value theory: Models and estimation. *Biometrika* **75** 397–415. MR0967580

Tilburg University
Department of Econometrics
and OR and CentER
P.O. Box 90153
NL-5000 LE Tilburg
The Netherlands
E-mail: j.h.j.einmahl@uvt.nl

Université catholique de Louvain
Institut de statistique
Voie du Roman Pays, 20
B-1348 Louvain-la-Neuve
Belgium
E-mail: johan.segers@uclouvain.be